\documentclass[12pt]{article}
 
\newcommand{\ignore}[1]{}

\usepackage{amsmath, amsthm, amsfonts, amssymb} 

\usepackage[dvips,final]{graphicx}

\usepackage{graphicx}

\usepackage{float}
\newfloat{figure}{H}{lof}
\floatname{figure}{\figurename}

\DeclareMathAlphabet{\eufrak}{U}{}{}{}  
\SetMathAlphabet\eufrak{normal}{U}{euf}{m}{n}
\SetMathAlphabet\eufrak{bold}{U}{euf}{b}{n}

\usepackage{amsmath, amsthm, amsfonts, amssymb}

\numberwithin{equation}{section}

\def\real{{\mathord{{\rm I\kern-2.8pt R}}}}        
\def\inte{{\mathord{{\rm I\kern-2.8pt N}}}}
\def\PP{{\mathord{{\rm I\kern-2.8pt P}}}}

\def\real{{\mathord{\mathbb R}}}

\def\inte{{\mathord{\mathbb N}}}

\def\R{\right}
\def\L{\left}

\newcommand{\disp}{\displaystyle}

\def\P{\mathbb{P}}
\def\E{\mathop{\hbox{\rm I\kern-0.20em E}}\nolimits}

\newtheorem{prop}{Proposition}[section]

\newtheorem{lemma}[prop]{Lemma}
\newtheorem{definition}[prop]{Definition}

\newtheorem{theorem}[prop]{Theorem}

\newtheorem{remarks}[prop]{Remarks}

\title{ 
\Huge Convergence of finite-dimensional laws of the weighted quadratic variations process for some fractional Brownian sheets
} 
 
\author{
\Large 
Anthony R\'eveillac\footnote{anthony.reveillac@univ-lr.fr}  
\\ 
Laboratoire Math\'ematiques Image et Applications
\\ 
Universit\'e de La Rochelle 
\\ 
Avenue Michel Cr\'epeau 
\\ 
17042 La Rochelle Cedex 
\\ 
France 
}

\allowdisplaybreaks

\begin{document}
\hyphenation{func-tio-nals} 
\hyphenation{para-me-ter}
\maketitle
\begin{abstract}
\noindent 
In this paper we state and prove a central limit theorem for the finite-dimensional laws of the quadratic variations process of certain fractional Brownian sheets. The main tool of this article is a method developed by Nourdin and Nualart in \cite{NourdinNualart} based on the Malliavin calculus.                                            
\end{abstract} 
 
\normalsize

\vspace{0.5cm}



\normalsize

\baselineskip0.7cm

\section{Introduction}
\label{section:introduction}

Recently in several works, the asymptotic behavior of the weighted $p$-power variations of stochastic processes has been investigated. For a one-parameter stochastic process $(Z_t)_{t\in[0,1]}$ observed at times $\{i/n,\; 0\leq i \leq n\}$ it is defined as
\begin{equation}
\label{eq:weightedppowervariations}
\sum_{i=1}^n f\L(Z_{\frac{i-1}{n}}\R) \L(Z_{\frac{i}{n}}-Z_{\frac{i-1}{n}}\R)^p,
\end{equation}
and one is interested by the convergence of this quantity as $n$ goes to infinity and by the nature of this convergence. Power variations play an important role in the recovering of statistical properties of discretely observed stochastic processes and as an example, the results obtained in this area have been applied in some financial econometric problems (see \cite{Barndorff-NielsenGraversenJacodShephard}) and in mathematical finance (\textit{e.g.} in \cite{AitSahaliaJacodEstimators,Jacod1,Jacod2}). In these settings the considered stochastic processes were one-parameter It\^
o semimartingales. \\\\
\noindent
In this paper we consider a different situation since we consider two-parameter processes which are not semimartingales. Actually, we state and prove a central limit theorem (Theorem \ref{theorem:mainresult}) for the finite-dimensional laws of the weighted quadratic variations process of certain fractional Brownian sheets. More precisely, let $(W_{(s,t)}^{\alpha,\beta})_{(s,t)\in[0,1]^2}$ be a fractional Brownian sheet with Hurst indices $\alpha$ and $\beta$ such that $0<\alpha<\frac12$, $0<\beta<\frac12$ with $\alpha+\beta>\frac12$ then for a weight function $f:\real \to \real$ satisfying hypothesis \textbf{(H)} presented below, we have,
\begin{equation}
\label{eq:res}
n^{-1} \sum_{i=1}^{[n \cdot]} \sum_{j=1}^{[n \bullet]} f\L(W_{\L(\frac{i-1}{n},\frac{j-1}{n}\R)}^{\alpha,\beta}\R) (n^{2(\alpha+\beta)} \vert \Delta_{i,j}W^{\alpha,\beta} \vert^2-1) \overset{fdd}{\underset{n\to\infty}{\longrightarrow}} \sigma_{\alpha,\beta} \int_0^\cdot \int_0^\bullet f\L(W_{(u,v)}^{\alpha,\beta}\R) dW_{(u,v)}
\end{equation} 
where $W$ is a standard Brownian sheet independent of $W^{\alpha,\beta}$ and the notation $fdd$ means that the convergence is in the sense of stable convergence of the finite-dimensional laws (see (\ref{eq:stabledefinition})). Note that the constant $\sigma_{\alpha,\beta}$ appearing in the limiting process can be expressed explicitly with respect to $\alpha$ and $\beta$.\\\\
\noindent
The study of $p$-power variations for fractional Brownian motion has been initiated by Gradinaru and Nourdin in \cite{GradinaruNourdin}, Neuenkirch and Nourdin in \cite{NeuenkirchNourdin} and by Nourdin in \cite{NourdinSDE2007} in view of obtaining exact rate of convergence of some approximating scheme of scalar stochastic differential equations driven by a fractional Brownian motion. This study has been recently pursued by Nourdin in \cite{NourdinQuadrCubic}, by Nourdin and Nualart in \cite{NourdinNualart} and by Nourdin, Nualart and Tudor in \cite{NourdinNualartTudor} where more references about this topic can be found. A complete description of the nature of the convergence of weighted $p$-power variations of the form (\ref{eq:weightedppowervariations}) for a fractional Brownian motion $B$ with Hurst index $H$ is given in \cite[Theorem 1]{NourdinNualartTudor}. More precisely central and non-central limit theorems are derived depending on the values of $p$ and $H$. Concerning the particular case of weighted quadratic variations ($p=2$) it is shown in \cite[Theorem 1]{NourdinNualartTudor} that for $\frac14<H<\frac34$,
$$ n^{-1/2} \sum_{i=1}^n f\L(B_{\frac{i-1}{n}}\R) \L(n^{2H}\vert B_{\frac{i}{n}}-B_{\frac{i-1}{n}} \vert^2-1\R) \overset{\mathcal{L}}{\underset{n\to\infty}{\longrightarrow}} \sigma_H \int_0^1 f(B_s) dW_s,$$
where $\sigma_H$ is an explicit constant depending on $H$ and $W$ is a standard Brownian motion independent of $B$. Note this result was also obtained in \cite{NourdinNualart} when $\frac14<H<\frac12$. For the two-parameter case a central limit theorem has been obtained in \cite{Reveillac} for the weighted quadratic variations process of a standard Brownian sheet and applied to the construction of an asymptotically normal estimator of the quadratic variation of a two-parameter diffusion process. However, in \cite{Reveillac}, the stable finite-dimensional convergence in law is obtained by using a result based on some semimartingales techniques developed by Jacod and Shiryaev in \cite{JacodShiryaev}. Consequently this method is useless for fractional Brownian sheets, that's why we propose to replace it by adapting an argument presented in \cite{NourdinNualart} based on the Malliavin calculus which is valid in general Gaussian context. In order to make this paper self-contained we have adapted and reproduced the main computations collected in Section \ref{section:useful lemmas} originally realized by Nourdin and Nualart in \cite{NourdinNualart}.\\\\
\noindent       
We proceed as follows. In Section \ref{section:preliminariesandnotations} we recall some definitions and properties of the fractional Brownian sheet and some elements of the Malliavin calculus relative to this process. Then in Section \ref{section:mainresults} we state and prove the convergence of the finite-dimensional laws of the weighted quadratic variations of certain fractional Brownian sheets (Theorem \ref{theorem:mainresult}). Technical arguments used in Section \ref{section:mainresults} are collected in Section \ref{section:useful lemmas}.

\section{Preliminaries and notations}
\label{section:preliminariesandnotations}

In this section we recall the definition of the fractional Brownian sheet and we present some elements of Malliavin calculus.  

\subsection{The fractional Brownian sheet}

Several extensions of the fractional Brownian motion have been proposed in the literature as for example the \textit{fractional Brownian field} (\textit{c.f.} \cite{Lindstrom,BonamiEstrade}), the \textit{L\'evy's fractional Brownian field} (\textit{c.f.} \cite{CiesielskiKamont}) and the \textit{fractional Brownian sheet} (\textit{c.f.} \cite{Kamont,AyacheLegerPontier2002}) we consider in this paper. The definitions and properties of this section can be found in \cite{AlosMazetNualart2001,TudorViens2003}.

\begin{definition}[Fractional Brownian sheet]
A fractional Brownian sheet $(W_{(s,t)}^{\alpha,\beta})_{(s,t)\in[0,1]^2}$ with Hurst indices $(\alpha,\beta)\in (0,1)^2$ is a centered two-parameter Gaussian process equal to zero on the set 
$$\{(s,t) \in [0,1]^2, \; s=0 \textrm{ or } t=0 \} $$
whose covariance function is given by,
\begin{eqnarray*}
\disp{R^{\alpha,\beta}((s_1,t_1),(s_2,t_2))}&:=&\disp{\E\L[ W_{(s_1,t_1)}^{\alpha,\beta} W_{(s_2,t_2)}^{\alpha,\beta} \R]}\\
&=&\disp{K^\alpha(s_1,s_2) K^\beta(t_1,t_2)}\\
&=&\disp{\frac12 \big( s_1^{2\alpha}+s_2^{2\alpha}-\vert s_1-s_2 \vert^{2 \alpha} \big) \frac12 \big( t_1^{2 \beta}+t_2^{2 \beta}-\vert t_1-t_2 \vert^{2 \beta} \big).}
\end{eqnarray*}
\end{definition}   
\noindent
We assume that $(W_{(s,t)}^{\alpha,\beta})_{(s,t)\in[0,1]^2}$ is defined on a complete probability space $\L(\Omega,\mathcal{F},\P\R)$ where $\mathcal{F}$ is generated by $W^{\alpha,\beta}$. Let us denote by $\mathcal{H}$ the reproducing kernel Hilbert space associated with $W^{\alpha,\beta}$, that is, $\mathcal{H}$ is the closure of the linear span generated by indicator functions on $[0,1]^2$ with respect to the scalar product
$$ \langle \textbf{1}_{[0,s_1]\times[0,t_1]}, \textbf{1}_{[0,s_2]\times[0,t_2]} \rangle_{\mathcal{H}} = R^{\alpha,\beta}((s_1,t_1),(s_2,t_2)).$$
The mapping $\disp{\textbf{1}_{[0,s]\times[0,t]}\mapsto W_{(s,t)}^{\alpha,\beta}}$ provides an isometry between $\mathcal{H}$ and the first chaos $H_1^{\alpha,\beta}$. For an element $\varphi$ of $\mathcal{H}$ we denote by $W^{\alpha,\beta}(\varphi)$ the image of $\varphi$ in $H_1^{\alpha,\beta}$.\\
Note that we can also give a representation of $(W_{(s,t)}^{\alpha,\beta})_{(s,t)\in[0,1]^2}$ as a stochastic integral of kernels $K^\alpha$ and $K^\beta$ with respect to a standard Brownian sheet $(W_{(s,t)})_{(s,t)\in[0,1]^2}$:
$$ W_{(s,t)}^{\alpha,\beta}=\int_0^s \int_0^t K^\alpha(s,u) K^\beta(t,v) \; dW_{(u,v)},\quad (s,t)\in[0,1]^2.$$
Using this representation, Tudor and Viens in \cite{TudorViens2003,TudorViens2006} have developed a Malliavin calculus with respect to $W^{\alpha,\beta}$. Note that in \cite{TudorViens2006} Tudor and Viens have also given an extension of the divergence integral. Now we present some elements of Malliavin calculus with respect to fractional Brownian sheets and especially the Malliavin integration by parts formula (\ref{eq:generalizedMallavinIBP}).

\subsection{Malliavin calculus for fractional Brownian sheet}

We recall some definitions and properties of the Malliavin calculus for the fractional Brownian sheet. These elements are contained in the general framework described in \cite{Nualart3} for Gaussian processes.\\\\
\noindent
For a cylindrical functional $F$ of the form
\begin{equation}
\label{eq:cylindric}
F=f\L(W^{\alpha,\beta}(\varphi_1),\ldots,W^{\alpha,\beta}(\varphi_n)\R), \quad n\geq 1, \; \varphi_1,\ldots,\varphi_n \in \mathcal{H}, \; f \in \mathcal{C}_b^\infty(\real^n),     
\end{equation}
we define the Malliavin derivative $DF$ of $F$ as,
$$ DF:=\sum_{i=1}^n \partial_i f\L(W^{\alpha,\beta}(\varphi_1),\ldots,W^{\alpha,\beta}(\varphi_n)\R) \varphi_i.$$
Furthermore $D:L^2(\Omega,\mathcal{F},\P) \to L^2(\Omega,\mathcal{F},\P;\mathcal{H})$ is a closable operator and it can be extended to the closure of Sobolev space $\mathbb{D}^{1,2}$ defined by functional $F$ whose norm $\|F\|_{1,2}$ is finite with,
$$\|F\|_{1,2}:=\E\L[ F^2 \R] + \E\L[\| DF \|_{\mathcal{H}}^2\R].$$
The adjoint operator $I_1$ of $D$ is defined by the following duality relationship 
$$ \E\L[ F I_1(u) \R] = \E\L[ \langle DF, u \rangle_{\mathcal{H}} \R],$$ 
for $F$ in $\mathbb{D}^{1,2}$ and for $u$ in $\mathcal{H}$ such that their exists $c_u>0$ verifying
$$ \L\vert \E\L[ \langle DG, u \rangle_{\mathcal{H}} \R] \R\vert \leq c_u \|G\|_{L^2(\Omega,\mathcal{F},\P)},\quad \textrm{ for every functional } G \textrm{ of the form } (\ref{eq:cylindric}).$$
Let $n\geq 1$. The $n$th Wiener chaos $\mathfrak{H}_n$ of $W^{\alpha,\beta}$ is the closed linear subspace of $L^2(\Omega,\mathcal{F},\P)$ generated by the random variables $\L\{H_n\L(W^{\alpha,\beta}(\varphi)\R), \; \varphi \in \mathcal{H}, \; \|\varphi\|_{\mathcal{H}}=1\R\}$ where $H_n$ denotes the $n$th Hermite polynomial. A linear isometry between the symmetric tensor product $\mathcal{H}^{\odot n}$ and $\mathfrak{H}_n$ is defined as,
\begin{equation}
\label{eq:I_nH_n}
I_n\L(\varphi^{\otimes n}\R) := n! H_n\L(W^{\alpha,\beta}(\varphi)\R).
\end{equation}
We conclude this section by the following integration by parts formula:
\begin{equation}
\label{eq:generalizedMallavinIBP}
\E\L[ F I_n(h) \R] = \E\L[ \langle D^n F, h \rangle_{\mathcal{H}^{\otimes n}} \R], \quad h\in \mathcal{H}^{\odot n}, \; F \in \mathbb{D}^{n,2},
\end{equation}
where $\mathbb{D}^{n,2}$ is the space of functionals $F$ such that $\|F\|_{n,2}$ is finite with
$$ \|F\|_{n,2}:=\E\L[ F^2 \R] + \sum_{i=1}^n \E\L[\| D^i F \|_{\mathcal{H}}^2\R]. $$

\section{Convergence of finite-dimensional laws}
\label{section:mainresults}

In this section we state and prove the central limit theorem (Theorem \ref{theorem:mainresult}) for the finite-dimensional laws of the weighted quadratic variations process of a fractional Brownian sheet $W^{\alpha,\beta}$ whose Hurst indices satisfy $0<\alpha<\frac12$, $0<\beta<\frac12$ with $\alpha+\beta>\frac12$. The main part of this result consist in adapting a method developed by Nourdin and Nualart in \cite{NourdinNualart} relying on the Malliavin calculus.\\\\
\noindent
Let $(W_{(s,t)}^{\alpha,\beta})_{(s,t)\in[0,1]^2}$ be a two-parameter fractional Brownian motion with Hurst indices $\alpha$ and $\beta$. Let $f:\real\to\real$ a function satisfying hypothesis \textbf{(H)} presented below. Let 
\begin{equation}
\label{eq:approximatingsequence}
X_{(s,t)}^n:=n^{-1} \sum_{i=1}^{[n s]} \sum_{j=1}^{[n t]} f\L(W_{\L(\frac{i-1}{n},\frac{j-1}{n}\R)}^{\alpha,\beta}\R) \L( n^{2(\alpha+\beta)} \vert\Delta_{i,j} W \vert^2 -1 \R),\quad (s,t)\in [0,1]^2
\end{equation}
be the re-normalized weighted quadratic variations where the increments $\Delta_{i,j}W^{\alpha,\beta}$ are defined as
$$\Delta_{i,j}W^{\alpha,\beta}:=W_{\L(\frac{i-1}{n},\frac{j-1}{n}\R)}^{\alpha,\beta}+W_{\L(\frac{i}{n},\frac{j}{n}\R)}^{\alpha,\beta}-W_{\L(\frac{i-1}{n},\frac{j}{n}\R)}^{\alpha,\beta}-W_{\L(\frac{i-1}{n},\frac{j}{n}\R)}^{\alpha,\beta}.$$   
Recall that using the definition of $I_2$ (\ref{eq:I_nH_n}) and the relation between Hermite polynomials (see for example \cite[(1.3) p.~5]{Nualart3}) $X^n$ can be expressed as,
$$X_{(s,t)}^n:=n^{-1} \sum_{i=1}^{[n s]} \sum_{j=1}^{[n t]} f\L(W_{\L(\frac{i-1}{n},\frac{j-1}{n}\R)}^{\alpha,\beta}\R) n^{2(\alpha+\beta)} I_2(\textbf{1}_{\Delta_{i,j}}^{\otimes 2}),\quad (s,t)\in [0,1]^2.$$
We aim at showing that $X^n$ converges stably in finite-dimensional law to a process $X$ defined as
\begin{equation}
\label{eq:limitingprocess}
X_{(s,t)}:=\sigma_{\alpha,\beta} \int_0^s \int_0^t f\L(W_{(u,v)}^{\alpha,\beta}\R) dB_{(u,v)},\quad (s,t)\in [0,1]^2,
\end{equation}
where $B$ is a standard Brownian sheet independent of $W^{\alpha,\beta}$ and
\begin{equation}
\label{eq:limitingconstant}
\sigma_{\alpha,\beta}:=\sqrt{ \frac18 \sum_{c,d=-\infty}^\infty \L(\L[\vert c+1\vert^{2\alpha} + \vert c-1\vert^{2\alpha} - 2 \vert c \vert^{2\alpha}\R] \L[\vert d+1\vert^{2 \beta} + \vert d-1\vert^{2 \beta} - 2 \vert d\vert^{2 \beta}\R] \R)^2 }.
\end{equation}
Note that in the standard Brownian sheet case ($(\alpha,\beta)=(\frac12,\frac12)$), $X$ coincide with the limiting process obtained in \cite{Reveillac} since $\sigma_{\frac12,\frac12}:=\sqrt{2}$. Note that in \cite{Reveillac} it is shown that the sequence $X^n$ is tight in the two-parameter Skorohod space which enabled us to obtain in the case $(\alpha,\beta)=(\frac12,\frac12)$ the convergence in law of the sequence $(X^n)_n$ to $X$ and not only the stable convergence of the finite-dimensional laws.\\
Stable convergence in law has been introduced by R\'enyi in \cite{Renyi58,Renyi63}. A full explanation about this subject can be found in \cite[Section VIII.5.c]{JacodShiryaev}. Roughly speaking stable convergence in law consists in the convergence of the couple $(X^n,W)_n$ to $(X,W)$ as $n$ goes to infinity. The definition of the limiting process $X$ in (\ref{eq:limitingprocess}) requires some justification. A way to define $X$ is to consider an extension of the probability basis say $\mathcal{B}:=(\Omega,\mathcal{F},(\mathcal{F}_{(s,t)})_{(s,t)\in [0,1]^2},\P)$ on which $W^{\alpha,\beta}$ is defined. We introduce an auxiliary probability basis $\mathcal{B}':=\mathcal{B}:=(\Omega',\mathcal{F}',(\mathcal{F}_{(s,t)}')_{(s,t)\in [0,1]^2},\P')$ as follows. The set $\Omega'$ is defined as the space of continuous functions on $[0,1]^2$ which vanish on $\{(s,t)\in [0,1]^2, \; s=0 \textrm{ or } t=0\}$. The process $B$ is then the canonical process on $\Omega'$ that is 
$$B_{(s,t)}(\omega'):=\omega'(s,t),\quad (s,t)\in[0,1]^2$$ 
and $\mathcal{F}'$ is defined as the filtration generated by $B$. According to \cite[Section 2.4.1]{Nualart3} there exists a probability measure $\P'$ on $(\Omega',\mathcal{F}')$ under which $B$ is a standard Brownian sheet. We finally define $\mathcal{F}_{(s,t)}'$ as the $\sigma$-field generated by $\{B_{(u,v)},\; (u,v)\preceq (s,t) \}$ with 
$$(s,t) \preceq (u,v) \Leftrightarrow s\leq u \textrm{ and } t\leq v.$$
We now describe the extension $\tilde{\mathcal{B}}:=(\tilde{\Omega},\tilde{\mathcal{F}},(\tilde{\mathcal{F}}_{(s,t)})_{(s,t) \in[0,1]^2},\tilde{\P})$ of $\mathcal{B}$ on which $X$ will be defined.
$$
\left\lbrace
\begin{array}{l}
\tilde{\Omega}:=\Omega \times \omega',\\
\tilde{\mathcal{F}}:=\mathcal{F} \otimes \mathcal{F}',\\
(\tilde{\mathcal{F}}_{(s,t)\in[0,1]^2}):=(\cap_{(u,v)>(s,t)} \mathcal{F}_{(u,v)} \otimes \mathcal{F}_{(s,t)}')_{(s,t) \in [0,1]^2},\\
\tilde{\P}(d\omega,d\omega'):=\P(d\omega) \tilde{\P}(d\omega').
\end{array}
\right.
$$ 
On $\tilde{\mathcal{B}}$, $(X_{(s,t)})_{(s,t)\in[0,1]^2}$ is defined by
$$ X_{(s,t)}(\omega,\omega'):=\sigma_{\alpha,\beta} \L( \int_0^s \int_0^t f\L(W_{(u,v)}^{\alpha,\beta}\R)(\omega) \; dB_{(u,v)} \R) (\omega'), \quad (s,t) \in [0,1]^2.$$
\noindent  
We will denote by $\E$ (respectively $\tilde{\E}$) the expectation under $\P$ (respectively $\tilde{\P}$).\\ 
Before stating and proving the first result we introduce the hypothesis (\textbf{H}) on the weight function $f$.\\
\textbf{Hypothesis (H)}:\\
$f:\real \to \real$ belongs to $\mathcal{C}^4(\real)$ and for every $i=1,\ldots,4$ 
$$\sup_{(s,t)\in[0,1]^2} \E\L[\L\vert f^{(i)}\L(W_{(s,t)}^{\alpha,\beta}\R )\R\vert^p \R]<\infty, \quad p\in(0,\infty) .$$

\begin{theorem}
\label{theorem:mainresult}
Let $f:\real\to\real$ satisfying hypothesis \emph{(\textbf{H})}. If one of the following situations hold
\begin{itemize}
\item[(i)] $\frac14<\alpha<\frac12$ and $\frac14<\beta<\frac12$,
\item[(ii)] $0<\alpha<\frac14$ and $\frac14<\beta<\frac12$ with $\alpha+\beta>\frac12$,
\item[(ii')] $0<\beta<\frac14$ and $\frac14<\alpha<\frac12$ with $\alpha+\beta>\frac12$.
\end{itemize}
then the sequence $(X^n)_n$ defined in (\ref{eq:approximatingsequence}) converges stably in finite-dimensional law to $X$ (see (\ref{eq:limitingprocess})), that is for every $(t_1,\ldots,t_m) \in ([0,1]^2)^m$, for every bounded measurable function $g:\real^m\to\real$ and for every $\mathcal{F}$-measurable random variable $Z$, the following convergence holds
\begin{equation}
\label{eq:stabledefinition}
\E\L[ g(X_{t_1}^n,\ldots,X_{t_m}^n) Z \R] \underset{n\to\infty}{\longrightarrow} \tilde{\E}\L[ g(X_{t_1},\ldots,X_{t_m}) Z \R].
\end{equation}
\end{theorem}
\noindent
Recall that $\mathcal{F}$ denotes the $\sigma$-field generated by $(W_{(s,t)}^{\alpha,\beta})_{(s,t)\in[0,1]^2}$.\\ 
Before proving this result, let us make the following remarks.
\begin{remarks}
\begin{itemize}
\item[1)] The set of conditions (i), (ii) and (ii') in Theorem \ref{theorem:mainresult} is equivalent to the condition
$$0 <\alpha< \frac12 \; \textrm{ and } \; 0<\beta<\frac12 \; \textrm{ with } \; \alpha+\beta>\frac12.$$
This last condition will be preferred in the computations of the proof.
\item[2)] Condition (ii) (and its symmetric counterpart (ii')) in Theorem \ref{theorem:mainresult} states that one of the two components of the sheet is allowed to have a regularity lower than $\frac14$ provided the other one has a regularity such that $\alpha+\beta>\frac12$. Observe that $\frac14$ is precisely the critical lower bound in the one-parameter case (see \cite[Theorem 1]{NourdinNualartTudor}). 
\item[3)] The standard Brownian case corresponding to $(\alpha,\beta)=(\frac12,\frac12)$ has been obtained in \cite[Theorem 3.1]{Reveillac} where in addition the tightness of $(X^n)_n$ in the two-parameter Skorohod space is proven. 
\end{itemize}
\end{remarks}
\noindent
Now we proof below Theorem \ref{theorem:mainresult}.
\begin{proof} 
Let $(\alpha,\beta)$ such that $0<\alpha<\frac12$, $0<\beta<\frac12$ with $\alpha+\beta>\frac12$. Let $m\geq 1$ and $t_1,\ldots,t_m$ in $([0,1]^2)^m$. For $1\leq i\leq m$, $t_i$ is denoted as $t_i:=(t_{i,1},t_{i,2})$. By Lemma \ref{lemma:tightness} the sequence $((X_{t_1}^n,\ldots,X_{t_m}^n),(W_z^{\alpha,\beta})_{z\in[0,1]^2})_n$ is tight in $\real^m \times \mathcal{C}([0,1]^2)$, in other words there exists a subsequence of $((X_{t_1}^n,\ldots,X_{t_m}^n),(W_z^{\alpha,\beta})_{z\in[0,1]^2})_n$ which converges in law to some limit denoted by $((X_{t_1}^\infty,\ldots,X_{t_m}^\infty),(W_z^{\alpha,\beta})_{z\in[0,1]^2})_n$. This subsequence and the sequence itself will be denoted by the same notation. We have to show that 
$$ \E\L[X_{t_1}^\infty,\ldots,X_{t_m}^\infty \vert (W_{z}^{\alpha,\beta})_{z\in[0,1]^2}\R] \overset{\mathcal{L}}{=} \E\L[X_{t_1},\ldots,X_{t_m}\vert (W_{z}^{\alpha,\beta})_{z\in[0,1]^2}\R].$$
Let us consider the conditional characteristic functions of respectively $(X_{t_1}^n,\ldots,X_{t_m}^n)$ and $(X_{t_1},\ldots,X_{t_m})$ given $W^{\alpha,\beta}$
\begin{equation}
\label{eq:definitioncaracteristicsfunctions}
\begin{array}{l}
\Phi^{m,n}{(\lambda_1,\ldots,\lambda_m)}:=\E\L[ e^{i \L\langle (\lambda_1,\ldots,\lambda_m),(X_{t_1}^n,\ldots,X_{t_m}^n)\R\rangle}\big\vert(W_z^{\alpha,\beta})_{z\in[0,1]^2} \R]\\
\Phi^m{(\lambda_1,\ldots,\lambda_m)}:=\E\L[ e^{i \L\langle (\lambda_1,\ldots,\lambda_m),(X_{t_1},\ldots,X_{t_m})\R\rangle} \big\vert(W_z^{\alpha,\beta})_{z\in[0,1]^2} \R].
\end{array}
\end{equation}
The proof consists to show that
$$ \Phi^{m,n}(\lambda_1,\ldots,\lambda_m) \underset{n\to\infty}{\longrightarrow} \Phi^{m}(\lambda_1,\ldots,\lambda_m), \quad \tilde{\P}-a.s.$$
Furthermore $ \Phi^{m}(\lambda_1,\ldots,\lambda_m)$ can be computed as follows,
$$  \Phi^{m}(\lambda_1,\ldots,\lambda_m)=e^{-\frac12 \langle (\lambda_1,\ldots,\lambda_m), \mathbb{Q}_{t_1,\ldots,t_m} (\lambda_1,\ldots,\lambda_m) \rangle_{\real^m}}, $$
where $\mathbb{Q}_{t_1,\ldots,t_m}$ is the symmetric matrix $\mathbb{Q}_{t_1,\ldots,t_m}:=(C_{i,j})_{1\leq i,j, \leq m}$ defined as,
$$ C_{i,j}:=\sigma_{\alpha,\beta}^2 \int_{[0,t_{i,1}\wedge t_{j,1}]\times[0,t_{i,2}\wedge t_{j,2}]}  f^2\L(W_{(u,v)}^{\alpha,\beta}\R) \; du dv,$$
which leads to,
$$ \Phi^m(\lambda_1,\ldots,\lambda_m)=\exp\L( -\frac12 \sum_{k=1}^m \lambda_k^2 C_{k,k} - \sum_{i=1}^{m-1} \lambda_i \sum_{k=i+1}^m \lambda_k C_{i,k} \R).$$
The key point is that $\Phi^m$ is the unique solution of the system of partial differential equations $(\ref{eq:s^m})$. Actually let the system of PDE's $(\ref{eq:s^m})$ given by functions $\varphi:\real^m \to \real$ satisfying,
\begin{equation}
\label{eq:s^m}
\left\lbrace
\begin{array}{ccc}
\frac{\partial \varphi}{\partial \lambda_1}(\lambda_1,\ldots,\lambda_m)&=&-\varphi(\lambda_1,\ldots,\lambda_m) \L( \lambda_1 C_{1,1} +\disp{\sum_{k=2}^m \lambda_k C_{1,k}} \R)\\
\vdots&\vdots&\vdots\\
\frac{\partial \varphi}{\partial \lambda_i}(\lambda_1,\ldots,\lambda_m)&=&-\varphi(\lambda_1,\ldots,\lambda_m) \L( \lambda_i C_{i,i} + \disp{\sum_{k=1; k\neq i}^m \lambda_k C_{i,k}} \R)\\
\vdots&\vdots&\vdots\\
\frac{\partial \varphi}{\partial \lambda_m}(\lambda_1,\ldots,\lambda_m)&=&-\varphi(\lambda_1,\ldots,\lambda_m) \L( \lambda_m C_{m,m} + \disp{\sum_{k=1}^{m-1} \lambda_k C_{m,k}} \R)
\end{array}
\right.,
\end{equation}
with the condition $\varphi(0,\ldots,0)=1$. Using standard techniques for PDE's systems that we briefly describe now, it can be shown that the unique solution to (\ref{eq:s^m}) is $\Phi^m$. Indeed, let $\varphi$  be a solution to (\ref{eq:s^m}). Using the first equation we have that
$$ \varphi(\lambda_1,\ldots,\lambda_m)=\exp\L(-\L(\frac{C_{11}}{2} \lambda_1^2 + \lambda_1 \sum_{k=2}^m C_{1,k} \lambda_k \R)\R) \exp\L(g_1(\lambda_2,\ldots,\lambda_m)\R),$$
where $g_1:\real^{m-1}\to\real$ is regular enough. Since $\varphi$ satisfies the second equation of the system we have,
$$ \frac{\partial g_1}{\partial \lambda_2}(\lambda_2,\ldots,\lambda_m)=-\L(\lambda_2 C_{2,2}+\sum_{k=3}^m \lambda_k C_{2,k}\R).$$
Consequently 
\begin{eqnarray*}
\disp{\varphi(\lambda_1,\ldots,\lambda_m)}&=&\disp{\exp\L(-\frac12 \sum_{i=1}^2 \lambda_i^2 C_{i,i} - \sum_{i=1}^2 \lambda_i \sum_{k=i+1}^m \lambda_k C_{i,k}\R)}\\
&&\disp{\times\exp\L(g_2(\lambda_3,\ldots,\lambda_m)\R)},
\end{eqnarray*}
where $g_2:\real^{m-2} \to \real$ is regular enough. Continuing this procedure we obtain that 
$$ \varphi(\lambda_1,\ldots,\lambda_m)=\exp\L(-\frac12 \sum_{i=1}^{m} \lambda_i^2 C_{i,i} - \sum_{i=1}^{m-1} \lambda_i \sum_{k=i+1}^m C_{i,k} \lambda_k \R) \exp(g_m),$$
where $g_m$ is a constant. Using the condition $\varphi(0,\ldots,0)=1$ we conclude that $g_m=0$ and consequently, $\varphi=\Phi^m$.\\
We proceed as follows: we show that the conditional characteristic function of $(X_{t_1}^\infty,\ldots,X_{t_m}^\infty)$ given $W^{\alpha,\beta}$ is solution to the system (\ref{eq:s^m}). This will be realized if we can show that for every random variable $H$ of the form $H:=\psi\L(W_{s_1}^{\alpha,\beta},\ldots,W_{s_r}^{\alpha,\beta}\R)$ with $\psi:\real^r \to \real$ in $\mathcal{C}_b^\infty(\real^r)$, $s_1,\ldots,s_r \in [0,1]^2$, we have
\begin{eqnarray}
\label{eq:interm2}
&&\disp{\frac{\partial}{\partial \lambda_p} \E\L[ e^{i \langle (\lambda_1,\ldots,\lambda_m), (X_{t_1}^\infty,\ldots,X_{t_m}^\infty) \rangle_{\real^m}} H \R]}\nonumber\\
&=&\disp{-\sigma_{\alpha,\beta}^2 \lambda_p \int_{[0,t_{p}]} \E\L[ f^2\L(W_{(s,t)}\R) H e^{i \langle \Lambda, \mathbb{X}^\infty \rangle} \R] ds dt}\nonumber\\ 
&-&\disp{\sigma_{\alpha,\beta}^2 \sum_{a=1,a\neq p}^m \lambda_a \int_{[0,t_{p,1}\wedge t_{a,1}]\times[0,t_{p,2}\wedge t_{a,2}]} \E\L[ f^2\L(W_{(s,t)}\R) H e^{i \langle \Lambda, \mathbb{X}^\infty \rangle} \R] ds dt,}
\end{eqnarray}
where we use the notation $\langle \Lambda,\mathbb{X}^{n} \rangle:=\sum_{k=1}^m \lambda_k X_{t_k}^{n}.$ Since
$$ \frac{\partial}{\partial \lambda_p} \E\L[ e^{i \langle (\lambda_1,\ldots,\lambda_m), (X_{t_1}^\infty,\ldots,X_{t_m}^\infty) \rangle_{\real^m}} H\R] = \lim_{n\to\infty} \frac{\partial \Phi_H^{m,n}}{\partial \lambda_p} (\lambda_1,\ldots,\lambda_m),\quad p=1,\ldots,m,$$
where $\Phi_H^{m,n} $ is defined as
$$\Phi_H^{m,n}:=\E\L[e^{\L(i \L\langle (\lambda_1,\ldots,\lambda_m), (X_{t_1}^n,\ldots,X_{t_m}^n) \R\rangle_{\real^m}\R)} H \R],$$
we have to compute $\frac{\partial \Phi_H^{m,n}}{\partial \lambda_p}(\lambda_1,\ldots,\lambda_m)$. The calculations presented below enable us to obtain the key expression (\ref{eq:finalexpressionofpartialderivative}).\\
Let $p\in \{1,\ldots,m\}.$ We have,
\begin{eqnarray}
\label{eq:expressionofpartialderivative}
&&\disp{ \frac{\partial \Phi_H^{m,n}}{\partial \lambda_p}(\lambda_1,\ldots,\lambda_m)}\nonumber\\
&=&\disp{i \E\L[ X_{t_p}^n e^{i \langle \Lambda,\mathbb{X}^{n} \rangle} H \R]}\nonumber\\
&=&\disp{ i n^{2(\alpha+\beta)-1} \sum_{i=1}^{[n t_{p,1}]} \sum_{j=1}^{[n t_{p,2}]} \E\L[ I_2(\textbf{1}_{\Delta_{i,j}}^{\otimes 2}) f\L(W_{\L(\frac{i-1}{n},\frac{j-1}{n}\R)}^{\alpha,\beta}\R) e^{i \langle \Lambda,\mathbb{X}^{n} \rangle} H \R] }\nonumber\\
&=& \disp{ i n^{2(\alpha+\beta)-1} \sum_{i=1}^{[n t_{p,1}]} \sum_{j=1}^{[n t_{p,2}]} \E\L[ \L\langle \textbf{1}_{\Delta_{i,j}}^{\otimes 2}, D^2\L(f\L(W_{\L(\frac{i-1}{n},\frac{j-1}{n}\R)}^{\alpha,\beta}\R) e^{i \langle \Lambda,\mathbb{X}^{n} \rangle} H\R) \R\rangle_{\mathcal{H}^{\otimes 2}} \R]}.
\end{eqnarray}
Let $1\leq i \leq [n t_{p,1}]$ and $1\leq j \leq [n t_{p,2}]$ and let $\delta_{i,j}:=\textbf{1}_{\L[0,\frac{i-1}{n}\R]\times\L[0,\frac{j-1}{n}\R]}$. We have,
\begin{eqnarray}
\label{eq:computaionexpectationD^2}
&&\disp{\E\L[ \L\langle \textbf{1}_{\Delta_{i,j}}^{\otimes 2}, D^2\L(f\L(W_{\L(\frac{i-1}{n},\frac{j-1}{n}\R)}^{\alpha,\beta}\R) e^{i \langle \Lambda,\mathbb{X}^n \rangle} H\R) \R\rangle_{\mathcal{H}^{\otimes 2}} \R]}\nonumber\\
&=&\disp{\E\L[ f''\L(W_{\L(\frac{i-1}{n},\frac{j-1}{n}\R)}^{\alpha,\beta}\R) H e^{i \langle \Lambda, \mathbb{X}^n \rangle} \R] \L\langle \delta_{i,j}, \textbf{1}_{\Delta_{i,j}} \R\rangle_{\mathcal{H}}^2 }\nonumber\\
&+&\disp{2 \E\L[ f'\L(W_{\L(\frac{i-1}{n},\frac{j-1}{n}\R)}^{\alpha,\beta}\R) e^{i \langle \Lambda, \mathbb{X}^n \rangle} \L\langle DH, \textbf{1}_{\Delta_{i,j}} \R\rangle_{\mathcal{H}} \R] \L\langle \delta_{i,j}, \textbf{1}_{\Delta_{i,j}} \R\rangle_{\mathcal{H}}}\nonumber\\
&+&\disp{ 2 i \E\L[ f'\L(W_{\L(\frac{i-1}{n},\frac{j-1}{n}\R)}^{\alpha,\beta}\R) H e^{i \langle \Lambda, \mathbb{X}^n \rangle} \L\langle D\langle \Lambda, \mathbb{X}^n\rangle, \textbf{1}_{\Delta_{i,j}} \R\rangle_{\mathcal{H}} \R] \L\langle \delta_{i,j}, \textbf{1}_{\Delta_{i,j}} \R\rangle_{\mathcal{H}} }\nonumber\\
&+&\disp{ \E\L[ f\L(W_{\L(\frac{i-1}{n},\frac{j-1}{n}\R)}^{\alpha,\beta}\R) e^{i \langle \Lambda, \mathbb{X}^n \rangle} \L\langle D^2H, \textbf{1}_{\Delta_{i,j}}^{\otimes 2} \R\rangle_{\mathcal{H}^\otimes 2} \R] }\nonumber\\
&+&\disp{ 2 i \E\L[ f\L(W_{\L(\frac{i-1}{n},\frac{j-1}{n}\R)}^{\alpha,\beta}\R) e^{i \langle \Lambda, \mathbb{X}^n \rangle} \L\langle D\L\langle \Lambda, \mathbb{X}^n\R\rangle, \textbf{1}_{\Delta_{i,j}} \R\rangle_{\mathcal{H}} \L\langle DH, \textbf{1}_{\Delta_{i,j}} \R\rangle_{\mathcal{H}} \R] }\nonumber\\
&+&\disp{ i \E\L[ f\L(W_{\L(\frac{i-1}{n},\frac{j-1}{n}\R)}^{\alpha,\beta}\R) H e^{i \langle \Lambda, \mathbb{X}^n \rangle} \L\langle D^2\langle \Lambda,\mathbb{X}^n \rangle, \textbf{1}_{\Delta_{i,j}}^{\otimes 2} \R\rangle_{\mathcal{H}^\otimes 2} \R] }\nonumber\\
&-&\disp{ \E\L[ f\L(W_{\L(\frac{i-1}{n},\frac{j-1}{n}\R)}^{\alpha,\beta}\R) H e^{i \langle \Lambda, \mathbb{X}^n \rangle} \L\langle D\langle \Lambda,\mathbb{X}^n \rangle, \textbf{1}_{\Delta_{i,j}} \R\rangle_{\mathcal{H}}^2 \R]. }
\end{eqnarray}
Now we compute $\L\langle D\langle \Lambda,\mathbb{X}^n \rangle, \textbf{1}_{\Delta_{i,j}} \R\rangle_{\mathcal{H}}$ and $\L\langle D^2\langle \Lambda,\mathbb{X}^n \rangle, \textbf{1}_{\Delta_{i,j}}^{\otimes 2} \R\rangle_{\mathcal{H}^{\otimes 2}}$.
\begin{eqnarray*}
&&\disp{\L\langle D\langle \Lambda,\mathbb{X}^n \rangle, \textbf{1}_{\Delta_{i,j}} \R\rangle_{\mathcal{H}}}\\
&=&\disp{ n^{2(\alpha+\beta)-1} \sum_{a=1}^m \lambda_a \sum_{k=1}^{[n t_{a,1}]} \sum_{l=1}^{[n t_{a,2}]} \bigg[ f'\L(W_{\L(\frac{k-1}{n},\frac{j-1}{n}\R)}^{\alpha,\beta}\R) I_2\L(\textbf{1}_{\Delta_{k,l}}^{\otimes 2}\R) \L\langle \delta_{k,l}, \textbf{1}_{\Delta_{i,j}} \R\rangle_{\mathcal{H}}}\\
&+&\disp{2 f\L(W_{\L(\frac{k-1}{n},\frac{l-1}{n}\R)}^{\alpha,\beta}\R) \Delta_{i,j} W^{\alpha,\beta} \L\langle \textbf{1}_{\Delta_{k,l}}, \textbf{1}_{\Delta_{i,j}} \R\rangle_{\mathcal{H}} \bigg] }
\end{eqnarray*}
\begin{eqnarray*}
&&\disp{\L\langle D^2\langle \Lambda,\mathbb{X}^n \rangle, \textbf{1}_{\Delta_{i,j}}^{\otimes 2} \R\rangle_{\mathcal{H}^\otimes 2}}\\
&=&\disp{ n^{2(\alpha+\beta)-1} \sum_{a=1}^m \lambda_a \sum_{k=1}^{[n t_{a,1}]} \sum_{l=1}^{[n t_{a,2}]} \bigg[ f''\L(W_{\L(\frac{k-1}{n},\frac{j-1}{n}\R)}^{\alpha,\beta}\R) I_2\L(\textbf{1}_{\Delta_{k,l}}^{\otimes 2}\R) \L\langle \delta_{k,l}, \textbf{1}_{\Delta_{i,j}} \R\rangle_{\mathcal{H}}^2 }\\
&+&\disp{ 4 f'\L(W_{\L(\frac{k-1}{n},\frac{l-1}{n}\R)}^{\alpha,\beta}\R) \Delta_{k,l} W^{\alpha,\beta} \L\langle \delta_{k,l}, \textbf{1}_{\Delta_{i,j}} \R\rangle_{\mathcal{H}}  \L\langle \textbf{1}_{\Delta_{k,l}}, \textbf{1}_{\Delta_{i,j}} \R\rangle_{\mathcal{H}} }\\
&+&\disp{ 2 f\L(W_{\L(\frac{k-1}{n},\frac{l-1}{n}\R)}^{\alpha,\beta}\R) \L\langle \textbf{1}_{\Delta_{k,l}}, \textbf{1}_{\Delta_{i,j}} \R\rangle_{\mathcal{H}}^2 \bigg] }
\end{eqnarray*}
In the following we use the following notation $[n t_a]=([n t_{a,1}],[n t_{a,2}])$. Equation (\ref{eq:computaionexpectationD^2}) can be rewritten as 
\begin{eqnarray}
\label{eq:expressionofpartialderivativebis}
&&\disp{\E\L[ \L\langle \textbf{1}_{\Delta_{i,j}}^{\otimes 2}, D^2\L(f\L(W_{\L(\frac{i-1}{n},\frac{j-1}{n}\R)}^{\alpha,\beta}\R) e^{i \langle \Lambda,\mathbb{X}^n \rangle} H\R) \R\rangle_{\mathcal{H}^{\otimes 2}} \R]}\nonumber\\
&=&\disp{2 i n^{2(\alpha+\beta)-1} \sum_{a=1}^m \lambda_a \sum_{k,l=1}^{[n t_{a}]} \E\L[ f\L(W_{\L(\frac{i-1}{n},\frac{j-1}{n}\R)}^{\alpha,\beta}\R) f\L(W_{\L(\frac{k-1}{n},\frac{l-1}{n}\R)}^{\alpha,\beta}\R) H e^{i \langle \Lambda, \mathbb{X}^n \rangle} \L\langle \textbf{1}_{\Delta_{k,l}}, \textbf{1}_{\Delta_{i,j}} \R\rangle_{\mathcal{H}}^2 \R] }\nonumber\\
&+&\disp{r_{i,j,n} },
\end{eqnarray}
where
\begin{eqnarray}
\label{eq:remainingterm}
\disp{r_{i,j,n}}&=&\disp{\E\L[ f''\L(W_{\L(\frac{i-1}{n},\frac{j-1}{n}\R)}^{\alpha,\beta}\R) H e^{i \langle \Lambda, \mathbb{X}^n \rangle} \R] \L\langle \delta_{i,j}, \textbf{1}_{\Delta_{i,j}} \R\rangle_{\mathcal{H}}^2}\nonumber\\
&+&\disp{2 \E\L[ f'\L(W_{\L(\frac{i-1}{n},\frac{j-1}{n}\R)}^{\alpha,\beta}\R) e^{i \langle \Lambda, \mathbb{X}^n \rangle} \L\langle DH, \textbf{1}_{\Delta_{i,j}} \R\rangle_{\mathcal{H}} \R] \L\langle \delta_{i,j}, \textbf{1}_{\Delta_{i,j}} \R\rangle_{\mathcal{H}}}\nonumber\\
&+&\disp{ 2 i \E\L[ f'\L(W_{\L(\frac{i-1}{n},\frac{j-1}{n}\R)}^{\alpha,\beta}\R) H e^{i \langle \Lambda, \mathbb{X}^n \rangle} \L\langle D\langle \Lambda, \mathbb{X}^n\rangle, \textbf{1}_{\Delta_{i,j}} \R\rangle_{\mathcal{H}} \R] \L\langle \delta_{i,j}, \textbf{1}_{\Delta_{i,j}} \R\rangle_{\mathcal{H}} }\nonumber\\
&+&\disp{ \E\L[ f\L(W_{\L(\frac{i-1}{n},\frac{j-1}{n}\R)}^{\alpha,\beta}\R) e^{i \langle \Lambda, \mathbb{X}^n \rangle} \L\langle D^2H, \textbf{1}_{\Delta_{i,j}}^{\otimes 2} \R\rangle_{\mathcal{H}^\otimes 2} \R] }\nonumber\\
&+&\disp{ 2 i \E\L[ f\L(W_{\L(\frac{i-1}{n},\frac{j-1}{n}\R)}^{\alpha,\beta}\R) e^{i \langle \Lambda, \mathbb{X}^n \rangle} \L\langle D\L\langle \Lambda, \mathbb{X}^n\R\rangle, \textbf{1}_{\Delta_{i,j}} \R\rangle_{\mathcal{H}} \L\langle DH, \textbf{1}_{\Delta_{i,j}} \R\rangle_{\mathcal{H}} \R] }\nonumber\\
&-&\disp{ \E\L[ f\L(W_{\L(\frac{i-1}{n},\frac{j-1}{n}\R)}^{\alpha,\beta}\R) H e^{i \langle \Lambda, \mathbb{X}^n \rangle} \L\langle D\langle \Lambda,\mathbb{X}^n \rangle, \textbf{1}_{\Delta_{i,j}} \R\rangle_{\mathcal{H}}^2 \R]}\nonumber\\
&+&\disp{i n^{2(\alpha+\beta)-1}\E\bigg[ f\L(W_{\L(\frac{i-1}{n},\frac{j-1}{n}\R)}^{\alpha,\beta}\R) H e^{i \langle \Lambda, \mathbb{X}^n \rangle} \sum_{a=1}^m \lambda_a \sum_{k,l=1}^{[n t_a]}}\nonumber\\ 
&&\disp{\times f''\L(W_{\L(\frac{k-1}{n},\frac{j-1}{n}\R)}^{\alpha,\beta}\R) I_2\L(\textbf{1}_{\Delta_{k,l}}^{\otimes 2}\R) \L\langle \delta_{k,l}, \textbf{1}_{\Delta_{i,j}} \R\rangle_{\mathcal{H}}^2 \bigg]}\nonumber\\ 
&+&\disp{4 i n^{2(\alpha+\beta)-1} \E\bigg[ f\L(W_{\L(\frac{i-1}{n},\frac{j-1}{n}\R)}^{\alpha,\beta}\R) H e^{i \langle \Lambda, \mathbb{X}^n \rangle} \sum_{a=1}^m \lambda_a} \nonumber\\
&&\disp{\times \sum_{k,l=1}^{[n t_a]} f'\L(W_{\L(\frac{k-1}{n},\frac{l-1}{n}\R)}^{\alpha,\beta}\R) \Delta_{k,l} W^{\alpha,\beta} \L\langle \delta_{k,l}, \textbf{1}_{\Delta_{i,j}} \R\rangle_{\mathcal{H}}  \L\langle \textbf{1}_{\Delta_{k,l}}, \textbf{1}_{\Delta_{i,j}} \R\rangle_{\mathcal{H}}\bigg]}\nonumber\\
&=&\disp{\sum_{k=1}^8 r_{i,j,n}^{(k)}.}
\end{eqnarray}
At this stage of the proof we have shown that,
\begin{eqnarray}
\label{eq:finalexpressionofpartialderivative}
&&\disp{\frac{\partial \Phi_H^{m,n}}{\partial \lambda_p} (\lambda_1,\ldots,\lambda_p)}\nonumber\\
&=&\disp{-2 n^{4(\alpha+\beta)-2} \sum_{a=1}^m \sum_{i,j=1}^{[n t_p]} \sum_{k,l=1}^{[n t_a]} \E\L[ f\L(W_{\L(\frac{i-1}{n},\frac{j-1}{n}\R)}^{\alpha,\beta}\R) f\L( W_{\L(\frac{k-1}{n},\frac{l-1}{n}\R)}^{\alpha,\beta}\R) H e^{i \langle \Lambda,\mathbb{X}^n \rangle } \L\langle \textbf{1}_{\Delta_{i,j}}, \textbf{1}_{\Delta_k,l} \R\rangle_{\mathcal{H}}^2 \R]}\nonumber\\
&+&\disp{ i n^{2(\alpha+\beta)-1} \sum_{i,j=1}^{[n t_p]} r_{i,j,n}.}
\end{eqnarray}
From Lemma \ref{lemma:estimates} we have
\begin{equation}
\label{eq:majorationoftheremainingterm}
\sup_{1\leq i,j\leq n}\vert r_{i,j,n} \vert \leq C n^{-4(\alpha+\beta)}, \quad n\geq 1.
\end{equation}
Using the estimate (\ref{eq:majorationoftheremainingterm}) we obtain,
\begin{footnotesize}
\begin{eqnarray*}
&&\disp{\lim_{n \to \infty} \frac{\partial \Phi_H^{m,n}}{\partial \lambda_p}(\lambda_1,\ldots,\lambda_m)}\\
&=&\disp{-2 \sum_{a=1}^m \lambda_a \lim_{n\to \infty} n^{4(\alpha+\beta)-2} \sum_{i,j=1}^{[n t_{p}]} \sum_{k,l=1}^{[n t_{a}]} \E\L[ f\L(W_{\L(\frac{i-1}{n},\frac{j-1}{n}\R)}^{\alpha,\beta}\R) f\L(W_{\L(\frac{k-1}{n},\frac{l-1}{n}\R)}^{\alpha,\beta}\R) H e^{i \langle \Lambda, \mathbb{X}^n \rangle} \L\langle \textbf{1}_{\Delta_{k,l}}, \textbf{1}_{\Delta_{i,j}} \R\rangle_{\mathcal{H}}^2 \R]}\\
&=&\disp{-\lim_{n\to \infty} \frac{1}{8 n^2} \sum_{a=1}^m \lambda_a \sum_{i,j=1}^{[n t_{p}]} \sum_{k,l=1}^{[n t_{a}]} \E\bigg[ f\L(W_{\L(\frac{i-1}{n},\frac{j-1}{n}\R)}^{\alpha,\beta}\R) f\L(W_{\L(\frac{k-1}{n},\frac{l-1}{n}\R)}^{\alpha,\beta}\R) H e^{i \langle \Lambda, \mathbb{X}^n \rangle} }\\
&&\disp{\times\L(\L[\vert k-i+1\vert^{2\alpha} + \vert k-i-1\vert^{2\alpha} - 2 \vert k-i\vert^{2\alpha}\R] \L[\vert l-j+1\vert^{2\alpha} + \vert l-j-1\vert^{2\alpha} - 2 \vert l-j\vert^{2\alpha}\R] \R)^2 \bigg]}\\
&=&\disp{ -\lim_{n\to \infty} \frac{1}{8 n^2} \lambda_p \sum_{i,j=1}^{[n t_{p}]} \sum_{k,l=1}^{[n t_{p}]} \E\bigg[ f\L(W_{\L(\frac{i-1}{n},\frac{j-1}{n}\R)}^{\alpha,\beta}\R) f\L(W_{\L(\frac{k-1}{n},\frac{l-1}{n}\R)}^{\alpha,\beta}\R) H e^{i \langle \Lambda, \mathbb{X}^n \rangle} }\\
&&\disp{\times\L(\L[\vert k-i+1\vert^{2\alpha} + \vert k-i-1\vert^{2\alpha} - 2 \vert k-i\vert^{2\alpha}\R] \L[\vert l-j+1\vert^{2\alpha} + \vert l-j-1\vert^{2\alpha} - 2 \vert l-j\vert^{2\alpha}\R] \R)^2 \bigg]}\\
&-&\disp{ \lim_{n\to \infty} \frac{1}{8 n^2} \sum_{a=1, a\neq p}^m \lambda_a \sum_{i,j=1}^{[n t_{p}]} \sum_{k,l=1}^{[n t_{a}]} \E\bigg[ f\L(W_{\L(\frac{i-1}{n},\frac{j-1}{n}\R)}^{\alpha,\beta}\R) f\L(W_{\L(\frac{k-1}{n},\frac{l-1}{n}\R)}^{\alpha,\beta}\R) H e^{i \langle \Lambda, \mathbb{X}^n \rangle} }\\
&&\disp{\times\L(\L[\vert k-i+1\vert^{2\alpha} + \vert k-i-1\vert^{2\alpha} - 2 \vert k-i\vert^{2\alpha}\R] \L[\vert l-j+1\vert^{2\alpha} + \vert l-j-1\vert^{2\alpha} - 2 \vert l-j\vert^{2\alpha}\R] \R)^2 \bigg]}\\
&=&\disp{ -\lim_{n\to \infty} \frac{1}{8 n^2} \lambda_p \sum_{c=-\infty}^\infty \sum_{i=1\vee(1-c)}^{([n t_{p,1}]-c)\wedge[n t_{p,1}]} \sum_{d=-\infty}^\infty \sum_{j=1\vee(1-d)}^{([n t_{p,2}]-d)\wedge[n t_{p,2}]} \E\bigg[ f\L(W_{\L(\frac{i-1}{n},\frac{j-1}{n}\R)}^{\alpha,\beta}\R) f\L(W_{\L(\frac{i+c-1}{n},\frac{j+d-1}{n}\R)}^{\alpha,\beta}\R) }\\
&&\disp{\times H e^{i \langle \Lambda, \mathbb{X}^n \rangle} \L(\L[\vert c+1\vert^{2\alpha} + \vert c-1\vert^{2\alpha} - 2 \vert c \vert^{2\alpha}\R] \L[\vert d+1\vert^{2\alpha} + \vert d-1\vert^{2\alpha} - 2 \vert d\vert^{2\alpha}\R] \R)^2 \bigg]}\\
&-&\disp{ \lim_{n\to \infty} \frac{1}{8 n^2} \sum_{a=1, a\neq p}^m \lambda_a \sum_{c,d=-\infty}^\infty \sum_{i=1\vee(1-c)}^{([n t_{a,1}]-c)\wedge[n t_{p,1}]} \sum_{j=1\vee(1-d)}^{([n t_{a,2}]-d)\wedge[n t_{p,2}]}\E\bigg[ f\L(W_{\L(\frac{i-1}{n},\frac{j-1}{n}\R)}^{\alpha,\beta}\R) f\L(W_{\L(\frac{i+c-1}{n},\frac{j+d-1}{n}\R)}^{\alpha,\beta}\R)}\\
&&\disp{ \times H e^{i \langle \Lambda, \mathbb{X}^n \rangle} \L(\L[\vert c+1\vert^{2\alpha} + \vert c-1\vert^{2\alpha} - 2 \vert c \vert^{2\alpha}\R] \L[\vert d+1\vert^{2\alpha} + \vert d-1\vert^{2\alpha} - 2 \vert d \vert^{2\alpha}\R] \R)^2 \bigg]}\\
&=&\disp{-\sigma_{\alpha,\beta}^2 \lambda_p \int_{[0,t_{p}]} \E\L[ f^2\L(W_{(s,t)}\R) H e^{i \langle \Lambda, \mathbb{X}^\infty \rangle} \R] ds dt}\\ 
&-&\disp{\sigma_{\alpha,\beta}^2 \sum_{a=1,a\neq p}^m \lambda_a \int_{[0,t_{p,1}\wedge t_{a,1}]\times[0,t_{p,2}\wedge t_{a,2}]} \E\L[ f^2\L(W_{(s,t)}\R) H e^{i \langle \Lambda, \mathbb{X}^\infty \rangle} \R] ds dt,}\\
\end{eqnarray*}
\end{footnotesize}
which leads to (\ref{eq:interm2}).
As a consequence,
\begin{eqnarray*}
&&\disp{ \frac{\partial}{\partial_p} \E\L[ e^{i \langle \Lambda, \mathbb{X}^\infty \rangle} \vert (W_z^{\alpha,\beta})_{z\in[0,1]^2} \R] }\\
&=&\disp{-\sigma_{\alpha,\beta}^2 \lambda_p \int_{[0,t_{p}]} f^2\L(W_{(s,t)}\R) \E\L[ e^{i \langle \Lambda, \mathbb{X}^\infty \rangle} \vert (W_z^{\alpha,\beta})_{z\in[0,1]^2}\R] ds dt}\\ 
&-&\disp{\sigma_{\alpha,\beta}^2 \sum_{a=1,a\neq p}^m \lambda_a \int_{[0,t_{p,1}\wedge t_{a,1}]\times[0,t_{p,2}\wedge t_{a,2}]} f^2\L(W_{(s,t)}\R) \E\L[ e^{i \langle \Lambda, \mathbb{X}^\infty \rangle} \vert (W_z^{\alpha,\beta})_{z\in[0,1]^2}\R] ds dt} 
\end{eqnarray*}
which is exactly the $p$th equation of the system (\ref{eq:s^m}).
\end{proof}

\section{Useful lemmas}
\label{section:useful lemmas}

All the results of this section mimic those obtained by Nourdin and Nualart in \cite{NourdinNualart}. We stress that we adapt these computations to the two-parameter case in order to make this paper self-contained. Furthermore the conditions imposed on $\alpha$ and $\beta$ in Theorem \ref{theorem:mainresult} can be easily deduced from the calculations realized in this section. Note that in this section $C$ will denote a generic constant which can differ from one line to another.

\begin{lemma}
\label{lemma:maj}
Let $l=(l_1,l_2)\in [0,1]^2$, $s=(s_1,s_2) \preceq t=(t_1,t_2).$ If $\alpha<\frac12$ and $\beta<\frac12$ then,
$$\L\vert\E\L[ W_{(l_1,l_2)}^{\alpha,\beta} \L(W_{(s_1,s_2)}^{\alpha,\beta}+W_{(t_1,t_2)}^{\alpha,\beta}-W_{(s_1,t_2)}^{\alpha,\beta}-W_{t_1,s_2}^{\alpha,\beta} \R) \R]\R\vert\leq \vert t_1-s_1\vert^{2\alpha} \vert t_2-s_2\vert^{2\beta}.$$
\end{lemma}

\begin{proof}
We refer to \cite[Lemma 1]{NourdinNualart}(or to \cite[Lemmas 4 and 5]{NourdinNualartTudor}) where the following estimate is shown,
\textrm{for } $a<b$ we have,
$$ \vert b^{2\gamma}-a^{2\gamma} \vert \leq \vert b-a \vert^{2\gamma}, \quad \forall \; 0<\gamma<\frac12. $$
\end{proof}

\begin{lemma}
\label{lemma:tightness}
Let $m\geq1$ and $t_1,\ldots,t_m$ in $[0,1]^2$. Under hypotheses of Theorem \ref{theorem:mainresult}, the sequence $(X_{t_1}^n,\ldots,X_{t_m}^n)_{n\geq 1}$ is tight in $\real^m$.
\end{lemma}

\begin{proof}
We show that 
\begin{equation}
\label{eq:tightestimate1}
\lim_{n\to\infty} \E\L[ (X_{t_1}^n,\ldots,X_{t_m}^n) \R]=0,
\end{equation}
\begin{equation}
\label{eq:tightestimate2}
\lim_{n\to\infty} \E\L[ \|(X_{t_1}^n,\ldots,X_{t_m}^n)\|^2 \R]=\sigma_{\alpha,\beta}^2 \sum_{i=1}^m \int_{[0,t_m]} \E\L[f^2\L(W_{(s,t)}^{\alpha,\beta}\R)\R] \, ds dt.
\end{equation}
\begin{small}\textbf{\textit{Proof of (\ref{eq:tightestimate1})}}\end{small}:\\
Since, $\Delta_{i,j}W^{\alpha,\beta}=I_1(\textbf{1}_{\Delta_{i,j}})$, using the definition $I_n(h^{\otimes n}):=n! H_n(W^{\alpha,\beta}(h))$ and relation on Hermite polynomials \cite[(1.3) p. 5]{Nualart3} 
$$ n^{2(\alpha+\beta)} \vert\Delta_{i,j}W^{\alpha,\beta}\vert^2 -1=n^{2(\alpha+\beta)} I_2(\textbf{1}_{\Delta_{i,j}}^{\otimes 2}).$$
Consequently,
\begin{eqnarray*}
\disp{\E\left[ f\L(W_{\L(\frac{k-1}{n},\frac{l-1}{n}\R)}^{\alpha,\beta}\R) I_2(\textbf{1}_{\Delta_{k,l}}^{\otimes 2}) \right]}&=&\disp{\E\left[ \left\langle D^2\left(f\L(W_{\L(\frac{k-1}{n},\frac{l-1}{n}\R)}^{\alpha,\beta}\R)\right),\textbf{1}_{\Delta_{k,l}}^{\otimes 2} \right\rangle_{\mathcal{H}^{\otimes 2}} \right]}\\
&=&\disp{\E\left[ f''\left(W_{\left(\frac{k-1}{n},\frac{l-1}{n}\right)}^{\alpha,\beta} \right) \right] \left\langle \delta_{k,l}^{\otimes 2},\textbf{1}_{\Delta_{k,l}}^{\otimes 2} \right\rangle_{\mathcal{H}^{\otimes 2}} }
\end{eqnarray*}
where $\delta_{k,l}=\textbf{1}_{\L[0,\frac{k-1}{n}\R]\times \L[0,\frac{l-1}{n}\R]}$.\\
Let for $p$ in $\{1,\ldots,n\}$ $t_p=(t_{p,1},t_{p,2})$. We have that
$$\left\vert \E\L[X_{t_p}^n\R] \right\vert \leq n^{2(\alpha+\beta)-1} \sum_{k=1}^{[n t_{p,1}]} \sum_{l=1}^{[n t_{p,2}]} \left\vert \E\left[ f''\left(W_{\left(\frac{k-1}{n},\frac{l-1}{n}\right)}^{\alpha,\beta} \right) \right] \right\vert \L\langle \delta_{k,l},\textbf{1}_{\Delta_{k,l}}\R\rangle_\mathcal{H}^2.$$
By Lemma \ref{lemma:maj} we have that
$$\left\vert \E\L[X_{t_p}^n\R] \right\vert \leq C \, n^{1-2(\alpha+\beta)}.$$
\begin{small}\textbf{\textit{Proof of (\ref{eq:tightestimate2})}}\end{small}:\\
Let $p$ in $\{1,\ldots,n\}$. 
\begin{equation}
\label{eq:a}
\vert X_{t_p}^n\vert^2= n^{4(\alpha+\beta)-2} \sum_{i,k=1}^{[n t_{p,1}]} \sum_{j,l=1}^{[n t_{p,2}]} f\L(W_{\L(\frac{i-1}{n},\frac{j-1}{n}\R)}^{\alpha,\beta}\R) f\L(W_{\L(\frac{k-1}{n},\frac{l-1}{n}\R)}^{\alpha,\beta}\R) I_2\L(\textbf{1}_{\Delta_{i,j}}^{\otimes 2}\R) I_2\L(\textbf{1}_{\Delta_{k,l}}^{\otimes 2}\R)
\end{equation}
Using the product formula (\cite[Proposition 1.1.2]{Nualart3}) we have,
\begin{eqnarray}
\label{eq:b}
\disp{I_2\L(\textbf{1}_{\Delta_{i,j}}^{\otimes 2}\R) I_2\L(\textbf{1}_{\Delta_{k,l}}^{\otimes 2}\R)}&=&\disp{\sum_{r=0}^{2} r! \L(\begin{array}{c}2\\r\end{array}\R) \L(\begin{array}{c}2\\r\end{array}\R) I_{4-2r}\L(\textbf{1}_{\Delta_{i,j}}^{\otimes 2} \otimes_r \textbf{1}_{\Delta_{k,l}}^{\otimes 2} \R) }\nonumber\\
&=& \disp{ I_4\L( \textbf{1}_{\Delta_{i,j}}^{\otimes 2} \tilde{\otimes} \textbf{1}_{\Delta_{k,l}}^{\otimes 2} \R) + 4 I_2\L( \textbf{1}_{\Delta_{i,j}} \otimes \textbf{1}_{\Delta_{k,l}} \R) \langle \textbf{1}_{\Delta_{i,j}}, \textbf{1}_{\Delta_{k,l}} \rangle_{\mathcal{H}}}\nonumber\\
&+&\disp{2 \langle \textbf{1}_{\Delta_{i,j}}^{\otimes 2}, \textbf{1}_{\Delta_{k,l}}^{\otimes 2} \rangle_{\mathcal{H}^{\otimes 2}}}.
\end{eqnarray}
From relations (\ref{eq:b}) and (\ref{eq:a}) we obtain that,
\begin{eqnarray}
\label{eq:c}
&&\disp{\E\L[\vert X_{t_p}^n\vert^2\R]}\\
&=&\disp{n^{4(\alpha+\beta)-2} \sum_{i,k=1}^{[n t_{p,1}]} \sum_{j,l=1}^{[n t_{p,2}]} \E\L[ f\L(W_{\L(\frac{i-1}{n},\frac{j-1}{n}\R)}^{\alpha,\beta}\R) f\L(W_{\L(\frac{k-1}{n},\frac{l-1}{n}\R)}^{\alpha,\beta}\R) I_4\L( \textbf{1}_{\Delta_{i,j}}^{\otimes 2} \tilde{\otimes} \textbf{1}_{\Delta_{k,l}}^{\otimes 2} \R) \R]}\nonumber\\
&+&\disp{4 n^{4(\alpha+\beta)-2} \sum_{i,k=1}^{[n t_{\alpha,1}]} \sum_{j,l=1}^{[n t_{\alpha,2}]} \E\L[ f\L(W_{\L(\frac{i-1}{n},\frac{j-1}{n}\R)}^{\alpha,\beta}\R) f\L(W_{\L(\frac{k-1}{n},\frac{l-1}{n}\R)}^{\alpha,\beta}\R) I_2\L( \textbf{1}_{\Delta_{i,j}} \otimes \textbf{1}_{\Delta_{k,l}} \R) \langle \textbf{1}_{\Delta_{i,j}}, \textbf{1}_{\Delta_{k,l}} \rangle_{\mathcal{H}}\R]}\nonumber\\
&+&\disp{ 2  n^{4(\alpha+\beta)-2} \sum_{i,k=1}^{[n t_{\alpha,1}]} \sum_{j,l=1}^{[n t_{\alpha,2}]} \E\L[ f\L(W_{\L(\frac{i-1}{n},\frac{j-1}{n}\R)}^{\alpha,\beta}\R) f\L(W_{\L(\frac{k-1}{n},\frac{l-1}{n}\R)}^{\alpha,\beta}\R) \langle \textbf{1}_{\Delta_{i,j}}^{\otimes 2}, \textbf{1}_{\Delta_{k,l}}^{\otimes 2} \rangle_{\mathcal{H}^{\otimes 2}}\R]}\nonumber\\
&=&\disp{T_1+T_2+T_3}.
\end{eqnarray}
We end the proof by considering independently each of the three terms in the right hand of expression (\ref{eq:c}).\\
First note that using Malliavin integration by parts formula the first term can be written as,
$$n^{4(\alpha+\beta)-2} \sum_{i,k=1}^{[n t_{p,1}]} \sum_{j,l=1}^{[n t_{p,2}]} \E\L[\L\langle D^4\L( f\L(W_{\L(\frac{i-1}{n},\frac{j-1}{n}\R)}^{\alpha,\beta}\R) f\L(W_{\L(\frac{k-1}{n},\frac{l-1}{n}\R)}^{\alpha,\beta}\R) \R), \textbf{1}_{\Delta_{i,j}}^{\otimes 2} \tilde{\otimes} \textbf{1}_{\Delta_{k,l}}^{\otimes 2} \R\rangle_{\mathcal{H}^{\otimes4}}\R] $$
which from Lemma \ref{lemma:maj} shows that $\vert T_1 \vert \leq C n^{2-4(\alpha+\beta)}$.\\
Then we consider the term $T_2$ in (\ref{eq:c}). We also integrate by part in order to obtain,
\begin{eqnarray*}
\disp{T_2}&=&\disp{ 4 n^{4(\alpha+\beta)-2} \sum_{i,k=1}^{[n t_{p,1}]} \sum_{j,l=1}^{[n t_{p,2}]} \E\L[ \L\langle D^2\L( f\L(W_{\L(\frac{i-1}{n},\frac{j-1}{n}\R)}^{\alpha,\beta}\R) f\L(W_{\L(\frac{k-1}{n},\frac{l-1}{n}\R)}^{\alpha,\beta}\R) \R), \textbf{1}_{\Delta_{i,j}} \otimes \textbf{1}_{\Delta_{k,l}} \R\rangle_{\mathcal{H}^{\otimes 2}} \R] }\\
&\times& \disp{\langle \textbf{1}_{\Delta_{i,j}}, \textbf{1}_{\Delta_{k,l}} \rangle_{\mathcal{H}} }\\
\end{eqnarray*}
As a consequence we obtain that,
\begin{eqnarray*}
\disp{\vert T_2 \vert} &\leq & \disp{C n^{-2} \sum_{i,k=1}^{[n t_{p,1}]} \sum_{j,l=1}^{[n t_{p,2}]} \vert \langle \textbf{1}_{\Delta_{i,j}}, \textbf{1}_{\Delta_{k,l}} \rangle_{\mathcal{H}} \vert}\\
&=&\disp{C n^{-2-2(\alpha+\beta)} \sum_{i,k=1}^{[n t_{p,1}]} \sum_{j,l=1}^{[n t_{p,2}]} \bigg[ \L\vert \vert k-i+1 \vert^{2\alpha} +\vert k-i-1 \vert^{2\alpha} -2 \vert k-i \vert^{2\alpha} \R\vert }\\
&&\disp{ \times\L\vert \vert l-j+1 \vert^{2 \beta} +\vert l-j-1 \vert^{2 \beta} -2 \vert l-j \vert^{2 \beta} \R\vert \bigg] }\\
&\leq& \disp{ C n^{-2(\alpha+\beta)} \sum_{c=-\infty}^{\infty} \sum_{d=-\infty}^{\infty} \L\vert \vert c+1 \vert^{2\alpha} +\vert c-1 \vert^{2\alpha} -2 \vert c \vert^{2\alpha} \R\vert \L\vert d+1 \vert^{2\beta} +\vert d-1 \vert^{2 \beta} -2 \vert d \vert^{2 \beta} \R\vert }.
\end{eqnarray*}
The serie presented above converges since $\alpha<\frac12$ and $\beta<\frac12$.
Consider now $T_3$.
\begin{footnotesize}
\begin{eqnarray*}
\disp{T_3}&=&\disp{ \frac18  n^{-2} \sum_{i,k=1}^{[n t_{p,1}]} \sum_{j,l=1}^{[n t_{p,2}]} \bigg[ \E\L[ f\L(W_{\L(\frac{i-1}{n},\frac{j-1}{n}\R)}^{\alpha,\beta}\R) f\L(W_{\L(\frac{k-1}{n},\frac{l-1}{n}\R)}^{\alpha,\beta}\R)\R] }\\
&&\disp{ \times\L( \L\vert \vert k-i+1 \vert^{2\alpha} +\vert k-i-1 \vert^{2\alpha} -2 \vert k-i \vert^{2\alpha} \R\vert \L\vert \vert l-j+1 \vert^{2 \beta} +\vert l-j-1 \vert^{2 \beta} -2 \vert l-j \vert^{2 \beta} \R\vert \R)^2 \bigg] }\\
&=&\disp{ \frac18  n^{-2} \sum_{c=-\infty}^\infty \sum_{i=1\vee(1-c)}^{([n t_{p,1}]-c)\wedge[n t_{p,1}]} \sum_{d=-\infty}^\infty \sum_{j=1\vee(1-d)}^{([n t_{p,2}]-d)\wedge[n t_{p,2}]} \E\bigg[ f\L(W_{\L(\frac{i-1}{n},\frac{j-1}{n}\R)}^{\alpha,\beta}\R) f\L(W_{\L(\frac{i+c-1}{n},\frac{j+d-1}{n}\R)}^{\alpha,\beta}\R) }\\
&&\disp{\times \L(\L[\vert c+1\vert^{2\alpha} + \vert c-1\vert^{2\alpha} - 2 \vert c \vert^{2\alpha}\R] \L[\vert d+1\vert^{2 \beta} + \vert d-1\vert^{2 \beta} - 2 \vert d\vert^{2 \beta}\R] \R)^2 \bigg]}\\
&\underset{n\to\infty} {\longrightarrow}&\disp{\sigma_{\alpha,\beta}^2 \int_{[0,t_p]} \E\L[ f^2\L(W_{(s,t)}^{\alpha,\beta}\R) \R] ds dt}.
\end{eqnarray*}
\end{footnotesize}
Note the series appearing in the term $T_3$ converges since $\alpha<\frac34$ and $\beta<\frac34$.
\end{proof}

\begin{lemma}
\label{lemma:estimates}
Let the notations of Theorem \ref{theorem:mainresult} and of its proof prevail. Under hypotheses of Theorem \ref{theorem:mainresult}, we have
$$\sup_{1\leq i,j\leq n}\vert r_{i,j,n} \vert \leq C n^{-4(\alpha+\beta)}, \quad n\geq 1.$$
\end{lemma}

\begin{proof}
Recall that in (\ref{eq:remainingterm}) $r_{i,j,n}$ is decomposed into eight terms,
$$r_{i,j,n}=\sum_{k=1}^8 r_{i,j,n}^{(k)}.$$
We show each of them is less or equal to $C n^{-4(\alpha+\beta)}$. From the definition of $H$ we deduce the following expressions 
$$
\left\lbrace
\begin{array}{l}
DH=\sum_{l=1}^r \partial_l \psi(W_{s_1}^{\alpha,\beta},\ldots,W_{s_r}^{\alpha,\beta}) \textbf{1}_{[0,s_l]}\\
D^2H=\sum_{l_2=1}^r \sum_{l_1=1}^r \partial_{l_2} \partial_{l_1} \psi(W_{s_1}^{\alpha,\beta},\ldots,W_{s_r}^{\alpha,\beta}) \textbf{1}_{[0,s_{l_2}]}\otimes\textbf{1}_{[0,s_{l_1}]}.
\end{array}
\right.
$$
By Lemma \ref{lemma:maj}, for $k=1,2,4$,
$$\vert r_{i,j,n}^{(k)}\vert \leq C n^{-4(\alpha+\beta)}.$$
\textbf{Estimate of }$r_{i,j,n}^{(6)}$:\\
First note that 
$$\vert r_{i,j,n}^{(6)}\vert \leq C \E\L[ \L(1+f^2\L(W_{\L(\frac{i-1}{n},\frac{j-1}{n}\R)}^{\alpha,\beta}\R)\R) \L\langle D\langle \Lambda, \mathbb{X}^n \rangle,\textbf{1}_{\Delta_{i,j}} \R\rangle_{\mathcal{H}}^2  \R].$$
Let $F_{i,j,n}:=1+f^2\L(W_{\L(\frac{i-1}{n},\frac{j-1}{n}\R)}^{\alpha,\beta}\R)$.
\begin{eqnarray*}
&&\disp{\E\L[ F_{i,j,n} \L\langle D\langle \Lambda, \mathbb{X}^n \rangle,\textbf{1}_{\Delta_{i,j}} \R\rangle_{\mathcal{H}}^2  \R]}\\
&\leq&\disp{n^{4(\alpha+\beta)-2} \sum_{a=1}^m \sum_{b=1}^m \lambda_a \lambda_b \sum_{k,l=1}^{[n t_a]} \sum_{\tilde{k},\tilde{l}=1}^{[n t_b]} \bigg[ A_{i,j,n}^{(1)}  \langle \delta_{k,l}, \textbf{1}_{\Delta_{i,j}} \rangle_{\mathcal{H}} \langle \delta_{\tilde{k},\tilde{l}}, \textbf{1}_{\Delta_{i,j}} \rangle_{\mathcal{H}}}\\
&&\disp{+ A_{i,j,n}^{(2)} \langle \textbf{1}_{\Delta_{k,l}}, \textbf{1}_{\Delta_{i,j}} \rangle_{\mathcal{H}} \langle \textbf{1}_{\Delta_{\tilde{k},\tilde{l}}}, \textbf{1}_{\Delta_{i,j}} \rangle_{\mathcal{H}} \bigg] }
\end{eqnarray*}
where
$$
\left\lbrace
\begin{array}{l}
A_{i,j,n}^{(1)}= 2 \E\L[ F_{i,j,n} f'\L(W_{\L(\frac{k-1}{n},\frac{l-1}{n}\R)}^{\alpha,\beta}\R) f'\L(W_{\L(\frac{\tilde{k}-1}{n},\frac{\tilde{l}-1}{n}\R)}^{\alpha,\beta}\R) I_2\L(\textbf{1}_{\Delta_{k,l}}^{\otimes 2}\R) I_2\L(\textbf{1}_{\Delta_{\tilde{k},\tilde{l}}}^{\otimes 2}\R) \R],\\\\
A_{i,j,n}^{(2)}=8 \E\L[ F_{i,j,n} f\L(W_{\L(\frac{k-1}{n},\frac{l-1}{n}\R)}^{\alpha,\beta}\R) f\L(W_{\L(\frac{\tilde{k}-1}{n},\frac{\tilde{l}-1}{n}\R)}^{\alpha,\beta}\R) \Delta_{k,l}W^{\alpha,\beta} \Delta_{\tilde{k},\tilde{l}}W^{\alpha,\beta} \R]. 
\end{array}
\right.
$$
The term $n^{4(\alpha+\beta)-2} \sum_{a=1}^m \sum_{b=1}^m \lambda_a \lambda_b \sum_{k,l=1}^{[n t_a]} \sum_{\tilde{k},\tilde{l}=1}^{[n t_b]} A_{i,j,n}^{(1)}$ is very similar to the term $\E[\vert X_{t_p}^n \vert^2]$ computed in Lemma \ref{lemma:tightness}. Furthermore by Lemma \ref{lemma:maj} we have
$$\L\vert \langle \delta_{k,l}, \textbf{1}_{\Delta_{i,j}} \rangle_{\mathcal{H}} \langle \delta_{\tilde{k},\tilde{l}}, \textbf{1}_{\Delta_{i,j}} \rangle_{\mathcal{H}} \R\vert \leq C n^{-4(\alpha+\beta)}$$
which leads to,
$$\L\vert n^{4(\alpha+\beta)-2} \sum_{a=1}^m \sum_{b=1}^m \lambda_a \lambda_b \sum_{k,l=1}^{[n t_a]} \sum_{\tilde{k},\tilde{l}=1}^{[n t_b]} \bigg[ A_{i,j,n}^{(1)}  \langle \delta_{k,l}, \textbf{1}_{\Delta_{i,j}} \rangle_{\mathcal{H}} \langle \delta_{\tilde{k},\tilde{l}}, \textbf{1}_{\Delta_{i,j}} \rangle_{\mathcal{H}}\bigg] \R\vert \leq C n^{-4(\alpha+\beta)}.$$ 
We consider the term $A_{i,j,n}^{(2)}$.
\begin{footnotesize}
\begin{eqnarray*}
&&\disp{\L\vert\E\L[ F_{i,j,n} f\L(W_{\L(\frac{k-1}{n},\frac{l-1}{n}\R)}^{\alpha,\beta}\R) f\L(W_{\L(\frac{\tilde{k}-1}{n},\frac{\tilde{l}-1}{n}\R)}^{\alpha,\beta}\R) \Delta_{k,l}W^{\alpha,\beta} \Delta_{\tilde{k},\tilde{l}}W^{\alpha,\beta} \R]\R\vert}\\
&=&\disp{\L\vert\E\L[ F_{i,j,n} f\L(W_{\L(\frac{k-1}{n},\frac{l-1}{n}\R)}^{\alpha,\beta}\R) f\L(W_{\L(\frac{\tilde{k}-1}{n},\frac{\tilde{l}-1}{n}\R)}^{\alpha,\beta}\R) I_1\L(\textbf{1}_{\Delta_{k,l}}\R) I_1\L(\textbf{1}_{\Delta_{\tilde{k},\tilde{l}}}\R) \R]\R\vert}\\
&\leq&\disp{\L\vert \E\L[ F_{i,j,n} f\L(W_{\L(\frac{k-1}{n},\frac{l-1}{n}\R)}^{\alpha,\beta}\R) f\L(W_{\L(\frac{\tilde{k}-1}{n},\frac{\tilde{l}-1}{n}\R)}^{\alpha,\beta}\R) I_2\L(\textbf{1}_{\Delta_{k,l}} \otimes \textbf{1}_{\Delta_{\tilde{k},\tilde{l}}}\R) \R] \R\vert}\\
&&\disp{+\L\vert \E\L[ F_{i,j,n} f\L(W_{\L(\frac{k-1}{n},\frac{l-1}{n}\R)}^{\alpha,\beta}\R) f\L(W_{\L(\frac{\tilde{k}-1}{n},\frac{\tilde{l}-1}{n}\R)}^{\alpha,\beta}\R) \L\langle \textbf{1}_{\Delta_{k,l}}, \textbf{1}_{\Delta_{\tilde{k},\tilde{l}}}\R\rangle_{\mathcal{H}} \R]  \R\vert ,\textrm{by \cite[Proposition 1.1.2]{Nualart3}}}\\
&=&\disp{ \L\vert \E\L[\L\langle D^2\L( F_{i,j,n} f\L(W_{\L(\frac{k-1}{n},\frac{l-1}{n}\R)}^{\alpha,\beta}\R) f\L(W_{\L(\frac{\tilde{k}-1}{n},\frac{\tilde{l}-1}{n}\R)}^{\alpha,\beta}\R) \R), \textbf{1}_{\Delta_{k,l}} \otimes \textbf{1}_{\Delta_{\tilde{k},\tilde{l}}}\R \rangle_{\mathcal{H}^{\otimes 2}} \R] \R\vert }\\
&&\disp{+\L\vert \E\L[ F_{i,j,n} f\L(W_{\L(\frac{k-1}{n},\frac{l-1}{n}\R)}^{\alpha,\beta}\R) f\L(W_{\L(\frac{\tilde{k}-1}{n},\frac{\tilde{l}-1}{n}\R)}^{\alpha,\beta}\R) \L\langle \textbf{1}_{\Delta_{k,l}}, \textbf{1}_{\Delta_{\tilde{k},\tilde{l}}}\R\rangle_{\mathcal{H}} \R]  \R\vert}\\
&\leq&\disp{C n^{-2(\alpha+\beta)}.}
\end{eqnarray*}
\end{footnotesize}
We don't give the details of the computations since the term 
$$\E\L[\L\langle D^2\L( F_{i,j,n} f\L(W_{\L(\frac{k-1}{n},\frac{l-1}{n}\R)}^{\alpha,\beta}\R) f\L(W_{\L(\frac{\tilde{k}-1}{n},\frac{\tilde{l}-1}{n}\R)}^{\alpha,\beta}\R) \R), \textbf{1}_{\Delta_{k,l}} \otimes \textbf{1}_{\Delta_{\tilde{k},\tilde{l}}}\R \rangle_{\mathcal{H}^{\otimes 2}} \R]$$ 
can be written as a sum of 
\begin{footnotesize}$$\E\L[F_{i,j,n} f^{(a)}\L(W_{\L(\frac{k-1}{n},\frac{l-1}{n}\R)}^{\alpha,\beta}\R) f^{(b)}\L(W_{\L(\frac{\tilde{k}-1}{n},\frac{\tilde{l}-1}{n}\R)}^{\alpha,\beta}\R) f^{(c)}\L(W_{\L(\frac{i-1}{n},\frac{j-1}{n}\R)}^{\alpha,\beta}\R) \L\langle \delta_{{u_1,v_1}} \otimes \delta_{{u_2,v_2}}, \textbf{1}_{\Delta_{k,l}} \otimes \textbf{1}_{\Delta_{\tilde{k},\tilde{l}}} \R\rangle_{\mathcal{H}^{\otimes 2}} \R],$$\end{footnotesize}
where $a,b,c$ are integer between zero and two and $(u_i,v_i)$ are for $(k,l)$ or $(\tilde{k},\tilde{l})$, $i=1,2$. Consequently,
\begin{eqnarray*}
&&\disp{\L\vert n^{4(\alpha+\beta)-2} \sum_{a=1}^m \sum_{b=1}^m \lambda_a \lambda_b \sum_{k,l=1}^{[n t_a]} \sum_{\tilde{k},\tilde{l}=1}^{[n t_b]} A_{i,j,n}^{(2)} \langle \textbf{1}_{\Delta_{k,l}}, \textbf{1}_{\Delta_{i,j}} \rangle_{\mathcal{H}} \langle \textbf{1}_{\Delta_{\tilde{k},\tilde{l}}}, \textbf{1}_{\Delta_{i,j}} \rangle_{\mathcal{H}} \R\vert }\\
&\leq&\disp{n^{2(\alpha+\beta-1)} \sum_{a=1}^m \sum_{b=1}^m \lambda_a \lambda_b \sum_{k,l=1}^{[n t_a]} \sum_{\tilde{k},\tilde{l}=1}^{[n t_b]} A_{i,j,n}^{(2)} \L\vert \langle \textbf{1}_{\Delta_{k,l}}, \textbf{1}_{\Delta_{i,j}} \rangle_{\mathcal{H}} \langle \textbf{1}_{\Delta_{\tilde{k},\tilde{l}}}, \textbf{1}_{\Delta_{i,j}} \rangle_{\mathcal{H}} \R\vert}\\
&=&\disp{C n^{-2(\alpha+\beta)-2} \L(\sum_{c=-\infty}^\infty \vert \vert c+1 \vert^{2\alpha}+\vert c-1 \vert^{2\alpha} -2\vert c \vert^{2\alpha} \vert\R)^2}\\
&&\disp{\times \L(\sum_{d=-\infty}^\infty \vert \vert d+1 \vert^{2\beta}+\vert d-1 \vert^{2\beta} -2\vert d \vert^{2\beta} \vert\R)^2}\\
&\leq&\disp{ C n^{-4(\alpha+\beta)},}
\end{eqnarray*}
where the series above converge since $0<\alpha<\frac12$ and $0<\beta<\frac12$. \\\\
\textbf{Estimate of }$r_{i,j,n}^{(3)}$ \textbf{and} $r_{i,j,n}^{(5)}$:\\
Note that 
$$ \vert r_{i,j,n}^{(3)} \vert + \vert r_{i,j,n}^{(5)} \vert \leq C n^{-2(\alpha+\beta)} \L\| \L\langle D\langle \Lambda, \mathbb{X}^n\rangle, \textbf{1}_{\Delta_{i,j}} \R\rangle_{\mathcal{H}} \R\|_{L^2(\Omega,\mathcal{F},\P)},$$
which shows that for $k=3,5$,
$$\vert r_{i,j,n}^{(k)} \vert \leq C n^{-4(\alpha+\beta)},$$
using estimates obtained previously in this proof.\\\\
\textbf{Estimate of }$r_{i,j,n}^{(7)}$:
We have,
\begin{eqnarray*}
&&\disp{\vert r_{i,j,n}^{(7)} \vert}\\ 
&\leq& \disp{C n^{-2(\alpha+\beta)-1} \sum_{a=1}^m \lambda_a \sum_{k,l=1}^{[n t_a]} \L\vert \E\L[ f\L(W_{\L(\frac{i-1}{n},\frac{j-1}{n}\R)}^{\alpha,\beta}\R) H e^{i \langle \Lambda, \mathbb{X}^n \rangle} f''\L(W_{\L(\frac{k-1}{n},\frac{j-1}{n}\R)}^{\alpha,\beta}\R) I_2\L(\textbf{1}_{\Delta_{k,l}}^{\otimes 2}\R) \R] \R\vert.}
\end{eqnarray*}
We use the integration by parts formula, 
\begin{footnotesize}
\begin{eqnarray*}
&&\disp{\E\L[ f\L(W_{\L(\frac{i-1}{n},\frac{j-1}{n}\R)}^{\alpha,\beta}\R) H e^{i \langle \Lambda, \mathbb{X}^n \rangle} f''\L(W_{\L(\frac{k-1}{n},\frac{j-1}{n}\R)}^{\alpha,\beta}\R) I_2\L(\textbf{1}_{\Delta_{k,l}}^{\otimes 2}\R) \R]}\\
&=&\disp{\E\L[ \L\langle D^2 \L(f\L(W_{\L(\frac{i-1}{n},\frac{j-1}{n}\R)}^{\alpha,\beta}\R) H e^{i \langle \Lambda, \mathbb{X}^n \rangle} f''\L(W_{\L(\frac{k-1}{n},\frac{j-1}{n}\R)}^{\alpha,\beta}\R)\R), \textbf{1}_{\Delta_{k,l}}^{\otimes 2} \R\rangle_{\mathcal{H}^{\otimes 2}} \R]}\\
&=&\disp{ \E\L[ f''\L(W_{\L(\frac{i-1}{n},\frac{j-1}{n}\R)}^{\alpha,\beta}\R) H e^{i \langle \Lambda, \mathbb{X}^n \rangle} f''\L(W_{\L(\frac{k-1}{n},\frac{j-1}{n}\R)}^{\alpha,\beta}\R)\R] \L\langle \delta_{k,l} \otimes \delta_{i,j}, \textbf{1}_{\Delta_{k,l}}^{\otimes 2} \R\rangle_{\mathcal{H}^{\otimes 2}} }\\
&&\disp{+2 \E\L[ f'\L(W_{\L(\frac{i-1}{n},\frac{j-1}{n}\R)}^{\alpha,\beta}\R) e^{i \langle \Lambda, \mathbb{X}^n \rangle} f''\L(W_{\L(\frac{k-1}{n},\frac{j-1}{n}\R)}^{\alpha,\beta}\R) \L\langle DH \tilde{\otimes} \delta_{i,j}, \textbf{1}_{\Delta_{k,l}}^{\otimes 2} \R\rangle_{\mathcal{H}^{\otimes 2}}\R]}\\
&&\disp{+2 i \E\L[ f'\L(W_{\L(\frac{i-1}{n},\frac{j-1}{n}\R)}^{\alpha,\beta}\R) H e^{i \langle \Lambda, \mathbb{X}^n \rangle} f''\L(W_{\L(\frac{k-1}{n},\frac{j-1}{n}\R)}^{\alpha,\beta}\R) \L\langle D \langle \Lambda, \mathbb{X}^n \rangle \tilde{\otimes} \delta_{i,j}, \textbf{1}_{\Delta_{k,l}}^{\otimes 2} \R\rangle_{\mathcal{H}^{\otimes 2}}\R]}\\
&&\disp{+2 \E\L[ f'\L(W_{\L(\frac{i-1}{n},\frac{j-1}{n}\R)}^{\alpha,\beta}\R) H e^{i \langle \Lambda, \mathbb{X}^n \rangle} f^{(3)}\L(W_{\L(\frac{k-1}{n},\frac{j-1}{n}\R)}^{\alpha,\beta}\R)\R] \L\langle \delta_{k,l} \tilde{\otimes} \delta_{i,j}, \textbf{1}_{\Delta_{k,l}}^{\otimes 2} \R\rangle_{\mathcal{H}^{\otimes 2}}}\\
&&\disp{+ \E\L[ f\L(W_{\L(\frac{i-1}{n},\frac{j-1}{n}\R)}^{\alpha,\beta}\R) e^{i \langle \Lambda, \mathbb{X}^n \rangle} f''\L(W_{\L(\frac{k-1}{n},\frac{j-1}{n}\R)}^{\alpha,\beta}\R) \L\langle D^2 H, \textbf{1}_{\Delta_{k,l}}^{\otimes 2} \R\rangle_{\mathcal{H}^{\otimes 2}}\R]}\\
&&\disp{+2 i \E\L[ f\L(W_{\L(\frac{i-1}{n},\frac{j-1}{n}\R)}^{\alpha,\beta}\R) e^{i \langle \Lambda, \mathbb{X}^n \rangle} f''\L(W_{\L(\frac{k-1}{n},\frac{j-1}{n}\R)}^{\alpha,\beta}\R) \L\langle D \langle \Lambda, \mathbb{X}^n \rangle \tilde{\otimes} DH, \textbf{1}_{\Delta_{k,l}}^{\otimes 2} \R\rangle_{\mathcal{H}^{\otimes 2}}\R]}\\
&&\disp{+2 \E\L[ f\L(W_{\L(\frac{i-1}{n},\frac{j-1}{n}\R)}^{\alpha,\beta}\R) e^{i \langle \Lambda, \mathbb{X}^n \rangle} f^{(3)}\L(W_{\L(\frac{k-1}{n},\frac{j-1}{n}\R)}^{\alpha,\beta}\R) \L\langle \delta_{k,l} \tilde{\otimes} DH, \textbf{1}_{\Delta_{k,l}}^{\otimes 2} \R\rangle_{\mathcal{H}^{\otimes 2}}\R]}\\
&&\disp{-\E\L[ f\L(W_{\L(\frac{i-1}{n},\frac{j-1}{n}\R)}^{\alpha,\beta}\R) H e^{i \langle \Lambda, \mathbb{X}^n \rangle} f''\L(W_{\L(\frac{k-1}{n},\frac{j-1}{n}\R)}^{\alpha,\beta}\R) \L\langle D \langle \Lambda, \mathbb{X}^n \rangle \otimes D \langle \Lambda, \mathbb{X}^n \rangle, \textbf{1}_{\Delta_{k,l}}^{\otimes 2} \R\rangle_{\mathcal{H}^{\otimes 2}}\R]}\\
&&\disp{+ i \E\L[ f\L(W_{\L(\frac{i-1}{n},\frac{j-1}{n}\R)}^{\alpha,\beta}\R) H e^{i \langle \Lambda, \mathbb{X}^n \rangle} f''\L(W_{\L(\frac{k-1}{n},\frac{j-1}{n}\R)}^{\alpha,\beta}\R) \L\langle D^2 \langle \Lambda, \mathbb{X}^n \rangle, \textbf{1}_{\Delta_{k,l}}^{\otimes 2} \R\rangle_{\mathcal{H}^{\otimes 2}}\R]}\\
&&\disp{+2 i \E\L[ f\L(W_{\L(\frac{i-1}{n},\frac{j-1}{n}\R)}^{\alpha,\beta}\R) H e^{i \langle \Lambda, \mathbb{X}^n \rangle} f^{(3)}\L(W_{\L(\frac{k-1}{n},\frac{j-1}{n}\R)}^{\alpha,\beta}\R) \L\langle D \langle \Lambda, \mathbb{X}^n \rangle \tilde{\otimes} \delta_{k,l}, \textbf{1}_{\Delta_{k,l}}^{\otimes 2} \R\rangle_{\mathcal{H}^{\otimes 2}}\R]}\\
&&\disp{+ \E\L[ f\L(W_{\L(\frac{i-1}{n},\frac{j-1}{n}\R)}^{\alpha,\beta}\R) H e^{i \langle \Lambda, \mathbb{X}^n \rangle} f^{(4)}\L(W_{\L(\frac{k-1}{n},\frac{j-1}{n}\R)}^{\alpha,\beta}\R)\R] \L\langle \delta_{k,l}^{\otimes 2}, \textbf{1}_{\Delta_{k,l}}^{\otimes 2} \R\rangle_{\mathcal{H}^{\otimes 2}}.}
\end{eqnarray*}
\end{footnotesize}
By Lemma \ref{lemma:maj} and the estimates shown above in this proof, 
$$\L\vert \langle  \delta_{s,t}, \textbf{1}_{\Delta_{k,l}} \rangle_{\mathcal{H}} \R\vert \leq C n^{-2(\alpha+\beta)}, \quad (s,t)\in\{(k,l),(i,j)\},$$
$$\L\|\L\langle DH, \textbf{1}_{\Delta_{k,l}} \R\rangle_{\mathcal{H}}\R\|_{L^2(\Omega,\mathcal{F},\P)}\leq C n^{-2(\alpha+\beta)},$$
and
$$\L\|\L\langle D^2 H, \textbf{1}_{\Delta_{k,l}}^{\otimes 2} \R\rangle_{\mathcal{H}^{\otimes 2}}\R\|_{L^2(\Omega,\mathcal{F},\P)}\leq C n^{-4(\alpha+\beta)}.$$
Furthermore as done in the part \textbf{Estimate of $r_{i,j,n}^{(6)}$},
$$\E\L[ \L(1+\L(f\L(W_{\L(\frac{i-1}{n},\frac{j-1}{n}\R)}^{\alpha,\beta}\R) f''\L(W_{\L(\frac{k-1}{n},\frac{j-1}{n}\R)}^{\alpha,\beta}\R)\R)^2\R) \L\langle D \langle \Lambda, \mathbb{X}^n \rangle, \textbf{1}_{\Delta_{k,l}} \R\rangle_{\mathcal{H}} \R] \leq C n^{-2(\alpha+\beta)}.$$ 
Consequently,
\begin{footnotesize}
$$ \vert r_{i,j,n}^{(7)} \vert \leq C n^{-2(\alpha+\beta)-1} \sum_{a=1}^m \lambda_a \L[n^{-4(\alpha+\beta)+2} + \sum_{k,l=1}^{[n t_a]} \L\| \L\langle D^2 \langle \Lambda, \mathbb{X}^n \rangle, \textbf{1}_{\Delta_{k,l}}^{\otimes 2} \R\rangle_{\mathcal{H}^{\otimes 2}} \R\|_{L^2(\Omega,\mathcal{F},\P)}\R].$$
$$\E\L[ \L\langle D^2 \langle \Lambda, \mathbb{X}^n \rangle, \textbf{1}_{\Delta_{k,l}}^{\otimes 2} \R\rangle_{\mathcal{H}^{\otimes 2}}^2 \R] \leq n^{4(\alpha+\beta)-2} \sum_{a,b=1}^m \lambda_a \lambda_b \sum_{i,j=1}^{[n t_a]} \sum_{\tilde{i},\tilde{j}=1}^{[n t_b]} \E[J_1 + J_2 + J_3],$$
\end{footnotesize}
where 
\begin{footnotesize}
$$
\L\lbrace
\begin{array}{l}
J_1:=3 f''\L(W_{\L(\frac{i-1}{n},\frac{j-1}{n}\R)}^{\alpha,\beta}\R) f''\L(W_{\L(\frac{\tilde{i}-1}{n},\frac{\tilde{j}-1}{n}\R)}^{\alpha,\beta}\R) I_2\L(\textbf{1}_{\Delta_{i,j}}^{\otimes 2}\R) I_2\L(\textbf{1}_{\Delta_{\tilde{i},\tilde{j}}}^{\otimes 2}\R) \L\langle \delta_{i,j},\textbf{1}_{\Delta_{k,l}} \R\rangle_{\mathcal{H}}^2 \L\langle \delta_{\tilde{i},\tilde{j}},\textbf{1}_{\Delta_{k,l}} \R\rangle_{\mathcal{H}}^2,\\
J_2:=48 f'\L(W_{\L(\frac{i-1}{n},\frac{j-1}{n}\R)}^{\alpha,\beta}\R) f'\L(W_{\L(\frac{\tilde{i}-1}{n},\frac{\tilde{j}-1}{n}\R)}^{\alpha,\beta}\R) \Delta_{i,j}W \Delta_{\tilde{i},\tilde{j}}W \L\langle \textbf{1}_{\Delta_{i,j}},\textbf{1}_{\Delta_{k,l}} \R\rangle_{\mathcal{H}} \L\langle \delta_{i,j},\textbf{1}_{\Delta_{k,l}} \R\rangle_{\mathcal{H}} \\ 
\quad \quad \times \L\langle \textbf{1}_{\Delta_{\tilde{i},\tilde{j}}},\textbf{1}_{\Delta_{k,l}} \R\rangle_{\mathcal{H}} \L\langle \delta_{\tilde{i},\tilde{j}},\textbf{1}_{\Delta_{k,l}} \R\rangle_{\mathcal{H}},\\ 
J_3:=12 f'\L(W_{\L(\frac{i-1}{n},\frac{j-1}{n}\R)}^{\alpha,\beta}\R) f'\L(W_{\L(\frac{\tilde{i}-1}{n},\frac{\tilde{j}-1}{n}\R)}^{\alpha,\beta}\R) \L\langle \textbf{1}_{\Delta_{i,j}}, \textbf{1}_{\Delta_{k,l}} \R\rangle_{\mathcal{H}}^2 \L\langle \textbf{1}_{\Delta_{\tilde{i},\tilde{j}}}, \textbf{1}_{\Delta_{k,l}} \R\rangle_{\mathcal{H}}^2.
\end{array}
\R.
$$
\end{footnotesize}
From computations made in the proof of Lemma \ref{lemma:tightness} we have that 
$$ n^{4(\alpha+\beta)-2} \sum_{a,b=1}^m \lambda_a \lambda_b \sum_{i,j=1}^{[n t_a]} \sum_{\tilde{i},\tilde{j}=1}^{[n t_b]} \E[J_1] \leq C n^{-8(\alpha+\beta)}.$$
Using an estimate obtained for the term $A_{i,j,n}^{(2)}$ obtained above in this proof we have that,
\begin{footnotesize}
\begin{eqnarray*}
&&\disp{\L\vert n^{4(\alpha+\beta)-2} \sum_{a,b=1}^m \lambda_a \lambda_b \sum_{i,j=1}^{[n t_a]} \sum_{\tilde{i},\tilde{j}=1}^{[n t_b]} \E[J_2] \R\vert}\\
&\leq& \disp{C n^{-2(\alpha+\beta)-2} \sum_{i,j=1}^{[n t_a]} \sum_{\tilde{i},\tilde{j}=1}^{[n t_b]} \L\vert \L\langle \textbf{1}_{\Delta_{i,j}},\textbf{1}_{\Delta_{k,l}} \R\rangle_{\mathcal{H}} \R\vert \L\vert \L\langle \textbf{1}_{\Delta_{i,j}},\textbf{1}_{\Delta_{k,l}} \R\rangle_{\mathcal{H}} \R\vert}\\
&\leq&\disp{C n^{-2(\alpha+\beta)-2} \L( \sum_{i,j=1}^{n} \L\vert \L\langle \textbf{1}_{\Delta_{i,j}},\textbf{1}_{\Delta_{k,l}} \R\rangle_{\mathcal{H}} \R\vert \R)^2 }\\
&=& \disp{C n^{-6(\alpha+\beta)-2} \L( \sum_{c,d=-\infty}^{\infty} \L\vert \vert c+1\vert^{2\alpha} + \vert c-1\vert^{2\alpha} -2\vert c \vert^{2\alpha}  \R\vert \L\vert \vert d+1\vert^{2\beta} + \vert d-1\vert^{2\beta} -2\vert d \vert^{2\beta}  \R\vert \R)^2}.
\end{eqnarray*}
\end{footnotesize}
Note that the serie above is convergent since $0<\alpha<\frac12$ and $0<\beta<\frac12$. For the last term we have,
\begin{footnotesize}
\begin{eqnarray*}
&&\disp{n^{4(\alpha+\beta)-2} \sum_{a,b=1}^m \lambda_a \lambda_b \sum_{i,j=1}^{[n t_a]} \sum_{\tilde{i},\tilde{j}=1}^{[n t_b]} \E[J_3]}\\
&\leq&\disp{n^{4(\alpha+\beta)-2} \L( \sum_{i,j=1}^{n} \L\langle \textbf{1}_{\Delta_{i,j}},\textbf{1}_{\Delta_{k,l}} \R\rangle_{\mathcal{H}}^2 \R)^2}\\
&\leq&\disp{n^{-4(\alpha+\beta)-2} \L( \sum_{c,d=-\infty}^{\infty} \L\vert \vert c+1\vert^{2\alpha} + \vert c-1\vert^{2\alpha} -2\vert c \vert^{2\alpha}  \R\vert^2 \L\vert \vert d+1\vert^{2\beta} + \vert d-1\vert^{2\beta} -2\vert d \vert^{2\beta}  \R\vert^2 \R)^2.}
\end{eqnarray*}
\end{footnotesize}
which leads to the result.\\\\
\textbf{Estimate of }$r_{i,j,n}^{(8)}$:
\begin{eqnarray*}
&&\disp{\E\L[ f\L(W_{\L(\frac{i-1}{n},\frac{j-1}{n}\R)}^{\alpha,\beta}\R) f'\L(W_{\L(\frac{k-1}{n},\frac{l-1}{n}\R)}^{\alpha,\beta}\R) H e^{i \langle \Lambda, \mathbb{X}^n \rangle} \Delta_{k,l}W^{\alpha,\beta} \R]}\\
&=&\disp{\E\L[ f'\L(W_{\L(\frac{i-1}{n},\frac{j-1}{n}\R)}^{\alpha,\beta}\R) f'\L(W_{\L(\frac{k-1}{n},\frac{l-1}{n}\R)}^{\alpha,\beta}\R) H e^{i \langle \Lambda, \mathbb{X}^n \rangle} \R] \L\langle \delta_{i,j},\textbf{1}_{\Delta_{k,l}} \R\rangle_{\mathcal{H}}}\\
&+&\disp{ \E\L[ f\L(W_{\L(\frac{i-1}{n},\frac{j-1}{n}\R)}^{\alpha,\beta}\R) f''\L(W_{\L(\frac{k-1}{n},\frac{l-1}{n}\R)}^{\alpha,\beta}\R) H e^{i \langle \Lambda, \mathbb{X}^n \rangle} \R] \L\langle \delta_{k,j},\textbf{1}_{\Delta_{k,l}} \R\rangle_{\mathcal{H}} }\\
&+&\disp{ \E\L[ f\L(W_{\L(\frac{i-1}{n},\frac{j-1}{n}\R)}^{\alpha,\beta}\R) f'\L(W_{\L(\frac{k-1}{n},\frac{l-1}{n}\R)}^{\alpha,\beta}\R) e^{i \langle \Lambda, \mathbb{X}^n \rangle} \L\langle DH,\textbf{1}_{\Delta_{k,l}} \R\rangle_{\mathcal{H}} \R]}\\
&+&\disp{i \E\L[ f\L(W_{\L(\frac{i-1}{n},\frac{j-1}{n}\R)}^{\alpha,\beta}\R) f'\L(W_{\L(\frac{k-1}{n},\frac{l-1}{n}\R)}^{\alpha,\beta}\R) H e^{i \langle \Lambda, \mathbb{X}^n \rangle} \L\langle D\langle \Lambda, \mathbb{X}^n \rangle ,\textbf{1}_{\Delta_{k,l}} \R\rangle_{\mathcal{H}} \R]}
\end{eqnarray*}
From preceding computations we have that
$$\L\vert r_{i,j,n}^{(8)} \R\vert \leq C n^{-2(\alpha+\beta)-1} \sum_{a=1}^m \sum_{k,l=1}^{[n t_a]} \L\vert \L\langle \textbf{1}_{\Delta_{i,j}}, \textbf{1}_{\Delta_{k,l}} \R\rangle_{\mathcal{H}} \R\vert \leq C n^{-4(\alpha+\beta)-1}$$
which concludes the proof.
\end{proof}

\section*{Acknowledgement}
I thank Ivan Nourdin for helpful discussions.\\
We thank the financial support of the \textit{European Social Fund (ESF)}.

\def\polhk#1{\setbox0=\hbox{#1}{\ooalign{\hidewidth
\lower1.5ex\hbox{`}\hidewidth\crcr\unhbox0}}} \def\cprime{$'$}


\begin{thebibliography}{1}
\footnotesize

\bibitem{AitSahaliaJacodEstimators}
Y.~A\"{i}t Sahalia and J.~Jacod.
\newblock Volatility estimators for discretely sampled L\'evy processes.
\newblock {\em Ann. Statist.}, 35(1), 355--392, 2007.

\bibitem{AlosMazetNualart2001}
E.~Al\`os, O.~Mazet and D.~Nualart.
\newblock Stochastic calculus with respect to Gaussian processes.
\newblock {\em Ann. Probab.}, 29(2), 766--801, 2001.

\bibitem{AyacheLegerPontier2002}
A.~Ayache, S.~Leger and M.~Pontier.
\newblock Drap brownien fractionnaire. (French) [The fractional Brownian sheet]
\newblock {\em Potential Anal.}, 17(1), 31--43, 2002.

\bibitem{Barndorff-NielsenGraversenJacodShephard}
O.~Barndorff-Nielsen, S.~Graversen and J.~Jacod, N.~Shephard.
\newblock Limit theorems for bipower variation in financial econometrics.
\newblock {\em Econometric Theory}, 22(4), 677--719, 2006.
          
\bibitem{Billingsley}
P.~Billingsley.
\newblock Convergence of probability measures, Second edition.
\newblock {\em Wiley Series in Probability and Statistics: Probability and Statistics.}, John Wiley And Sons, Inc., New York, 1999.

\bibitem{BonamiEstrade}
A.~Bonami and A.~Estrade.
\newblock Anisotropic analysis of some Gaussian models. 
\newblock {\em J. Fourier Anal. Appl.}, 9(3), 215--236, 2003.

\bibitem{CairoliWalsh}
R.~Cairoli and J.~Walsh.
\newblock Stochastic integrals in the plane.
\newblock {\em Acta Math.}, 134:111--183, 1975.

\bibitem{CiesielskiKamont}
Z.~Ciesielski and A.~Kamont.
\newblock L\'evy's fractional Brownian random field and function spaces.
\newblock {\em Acta Sci. Math.}, 60:99--118, 1995.

\bibitem{GradinaruNourdin}
M.~Gradinaru and I.~Nourdin.
\newblock {\em Weighted power variations of fractional Brownian motion
and application to approximating schemes.}
\newblock{\em Preprint}, 2007.

\bibitem{Jacod1}
J.~Jacod.
\newblock Asymptotic properties of realized power variations and related functionals of semimartingales.
\newblock To appear in {\em Stoch. Processes and their Appl.}, 2007.

\bibitem{Jacod2}
J.~Jacod.
\newblock Statistics and high frequency data.
\newblock Lecture given at the SEMSTAT 2007 held in La Manga (Spain), May 2007.

\bibitem{JacodShiryaev}
J.~Jacod and A.N.~Shiryaev.
\newblock {\em Limit theorems for stochastic processes. Second edition} Number 288 of
  {\em Grundlehren der Mathematischen Wissenschaften [Fundamental Principles of Mathematical Sciences]}.
\newblock Springer-Verlag, Berlin, 2003.

\bibitem{Kamont}
A.~Kamont.
\newblock On the fractional anisotropic Wiener field.
\newblock {\em Anna . Probab. Math. Statist.}, 16(1):85--98, 1996.

\bibitem{Lindstrom}
T.~ Lindstr\o m.
\newblock Fractional Brownian fields as integrals of white noise.
\newblock {\em Bull. London Math. Soc.}, 25(1):83--88, 1993.

\bibitem{NeuenkirchNourdin}
A.~Neuenkirch and I.~Nourdin.
\newblock {\em Exact rate of convergence of some approximation schemes associated to SDE's driven by a fractional Brownian motion.}
\newblock{\em J. Theor. Probab.}, 20(4):871-899, 2007.
 
\bibitem{NourdinQuadrCubic}
I.~Nourdin.
\newblock {\em Asymptotic behavior of certain weighted quadratic and cubic variations of fractional Brownian motion. }
\newblock To appear in {\em Ann. Probab.}, 2007.
  
\bibitem{NourdinSDE2007}
I.~Nourdin.
\newblock {\em A simple theory for the study of SDE's driven by a fractional Brownian motion, in dimension one.}
\newblock To appear in {\em S\'eminaire de Probabilit\'es \textbf{XLI}}, 2007. 
      
\bibitem{NourdinNualart}
I.~Nourdin and D.~Nualart.
\newblock {\em Central limit theorems for multiple Skorohod integrals.}
\newblock{\em Preprint}, 2008.

\bibitem{NourdinNualartTudor}
I.~Nourdin, D.~Nualart and C.A.~Tudor
\newblock {\em Central and non-central limit theorems for weighted power variations of fractional Brownian motion.}
\newblock {\em Preprint}, 2007.

\bibitem{Nualart3}
D.~Nualart.
\newblock {\em The Malliavin calculus and related topics} of
  {\em Probability and Its Applications}.
\newblock Springer Verlag, Berlin, Second edition, 2006.

\bibitem{Renyi58}
A.~R\'enyi.
\newblock On mixing sequences of sets.
\newblock {\em Acta Math. Acad. Sci. Hung.}, 9, 215--228, 1958.

\bibitem{Renyi63}
A.~R\'enyi.
\newblock On stable sequences of events.
\newblock {\em Sankhya, Ser A}, 25, 293--302, 1963.

\bibitem{Reveillac}
A.~R\'eveillac.
\newblock {\em Estimation of quadratic variation for two-parameter diffusions.}
\newblock{\em Preprint}, 2007.

\bibitem{TudorViens2003}
C.~Tudor and F.~Viens.
\newblock It\^o formula and local time for the fractional Brownian sheet.
\newblock {\em Electron. J. Probab.}, 8, 2003.

\bibitem{TudorViens2006}
C.~Tudor and F.~Viens.
\newblock It\^o formula for the two-parameter fractional Brownian motion using the extended divergence operator.
\newblock {\em Stochastics}, 78(6), 443--462, 2006.

\end{thebibliography}
\end{document}